\documentclass[10pt]{amsart}

\usepackage{graphicx}
\usepackage{xcolor}
\usepackage{url}

\newcommand{\Z}{{\mathbb Z}}
\newcommand{\Q}{{\mathbb Q}}

\newtheorem {theorem}{Theorem}
\newtheorem {Cor}{Corollary}
\newtheorem {Pro}{Proposition}
\newtheorem {Le}{Lemma}

\theoremstyle{definition}

\newtheorem{rem}{Remark}
\begin{document}

\title[Ankeny-Brauer-Chowla recurrences]{Common values of linear recurrences related to Ankeny-Brauer-Chowla polynomials}

\author[Attila~Peth\H{o}]{Attila~Peth\H{o}}
\address{Department of Computer Science,\\ University of Debrecen,\\
H-4002 Debrecen, P.O. Box 400}
\email{petho.attila@unideb.hu}

\begin{abstract}
Let $P_u(X)=XQ(X)(X-u)+1$, where $Q(X)\in \mathbb{Z}[X]$ is monic, separabel, has constant term $\pm 1$, and degree $q\ge 1$. Denote $(A_n(u))$ the linear recursive sequence having characteristic polynomial $P_u(X)$ and initial terms $A_0,\ldots,A_{q+1}\in \mathbb{Z}$ such that at least one of them is non-zero. 
Then there exist effectively computable constants $n_0,u_0$ depending only on $Q(X)$ and $\max\{|A_0|,\ldots, |A_{q+1}|\}$ such that if $|A_n(u)| = |A_m(u)|$ holds for some integers $u,n,m$ with $u>u_0$ and $n\not=m$ then $\max\{ |n|,|m|\} <n_0. $
\end{abstract}

\date{\today}

\maketitle

\section{Introduction}

Let $P_u(X) = X Q(X) (X-u) +1$ with $Q(X)\in \mathbb{Z}[X]$ and with the non-negative integer parameter $u$. Such kind of polynomials seems to be investigated at the first time by Ankeny, Brauer and Chowla in 1956 \cite{ABC}. They considered only such $Q$'s, which have solely simple, integer roots. They proved that $P_u(X)$ is irreducible and the number fields generated by one of its roots have as large as possible class number for infinitely many $u$. To the honor of the authors I will call $P_u(X)$ the family of ABC polynomials. 

One can associate to any integer polynomial a Thue equation. Denote the degree of $Q(X)$ by $q$ then its associated Thue equation is 
\begin{equation} \label{e:Thue}
Y^{q+2} P_u(X/Y) = X(Y^q Q(X/Y)) (X-u Y) + Y^{q+2} = 1.
\end{equation}

If $Q(X)=\prod_{j=1}^q (X-a_j)$ with $a_1,\ldots,a_q$ pairwise distinct integers, then \eqref{e:Thue} has the solutions $(1,0), (0,1), (a_j,1), \; j=1,\ldots,q$ and $(0,-1)$ if $q$ is even. We call them {\it trivial}. Bombieri and Schmidt \cite{BSch} proved that a Thue equation of degree $n$ admits at most $O(n)$ solutions. Moreover they pointed out that Thue equations associated to ABC polynomials are with this respect up to a constant factor extremal, because they admit, at least $n+1$ solutions.

For $q=1,\, a_1=1 $ and with $-(u+1)$ instead of $u$, Mignotte and Tzanakis \cite{MT} proved that if $|u|>3.67 \cdot 10^{32},$ then \eqref{e:Thue} has only one non-trivial solution, namely $(x,y) = (-u-1,-1)$. Mignotte \cite{Mig2000} proved the same for $|u| \ge 5$ and settled the cases $u=0,1,2,3,4$ too.

For $q=2,\, a_1=1,\, a_2=-1,$ equation \eqref{e:Thue} was completely solved for $|a| > 10^{28}$ by Peth\H{o} \cite{Pe}. His result was extended to $|a| \ge 2$ by Mignotte, Peth\H{o} and Roth \cite{MPR}. In this case, there is again exactly one (parametrized) non-trivial solution, namely $(x,y) = (1, -a)$. Marzenta \cite{Marzenta} gave an essentially new proof of the theorem of Peth\H{o}. His proof is based on a method developed for the solution of parametrized families of norm form equations and presented Amorozo, Masser and Zannier in \cite{AMZ}.

Thomas \cite{Th} investigated equation \eqref{e:Thue} in the case, when the factor $X-u$ disappear, but the $a_i$ denote integer polynomials of distinct degrees in an integer variable $u,\ \, a_i = a_i(u) \in \Z[u],\; i=1,2.$ He conjectured that \eqref{e:Thue} has only the trivial solutions if $u$ is large enough, and pointed out that the conjecture might also be true in specific cases of polynomials of equal degrees. In the case $n=3,\, a_1=u^k,\, a_2=u^l,$ where $0<k<l$ are integers, Thomas proved in \cite{Th}, that for $u \ge \bigl( 2 \cdot 10^6 \, (k+2l)\bigr)^{4.85/(l-k)},$ \eqref{e:Thue} has only the trivial solutions. Heuberger \cite{Heuberger} settled Thomas conjecture even if $Q(X)$ is the product of any number of factors $X - a_i(u)$, but only under certain assumptions.

The above examples show that not only the trivial, but most of the non-trivial solutions are polynomials in $u$. Zigler \cite{Ziegler} reformulated Thomas' conjecture for solutions  $(x,y)\in \mathbb{C}[u]$ of \eqref{e:Thue} and proved it even in the case when  $a_1(u), a_2(u) \in \mathbb{C}[u]$. 

The proofs most of the mentioned results consist of three main steps; the first is the reformulation of \eqref{e:Thue} to a weighted unit equation in two variables. In the second step it is shown that writing a solution of the weighted unit equation in an appropriate system of fundamental units, one of the exponents must be larger than $c u \log u$ with a positive constant $c$. In the final step some result on logarithmic forms imply that there is no solution provided $u$ is large enough. Halter-Koch, Lettl, Peth\H{o} and Tichy \cite{HLPT} implemented this program, when $P_u(X)$ denotes a classical ABC polynomial, i.e. when $Q(X)=\prod_{j=1}^q (X-a_j)$ with $a_1,\ldots,a_q$ pairwise distinct integers. They finished the first two steps successfully, but failed with the third one, because the available lower bounds for linear forms in logarithms of algebraic numbers were not sharp enough. They were able to finish the work by invoking the Lang-Waldschmidt conjecture, see \cite{La1}. Despite of the 25 years development the Lang-Waldschmidt conjecture is still far from a proof. 

The present paper was motivated by \cite{HLPT} and by the paper of Peth\H{o} and Tengely \cite{PT}. They investigated common values of linear recursive sequences associated to Shanks' simplest cubics. The Thue equation associated to the same family of polynomials was solved by Thomas \cite{Thomas} and Mignotte \cite{Mignotte1993} more than thirty years ago. Taking heart I tried the same game with ABC polynomials. More precisely write $P_u(X) = \sum_{j=0}^{q+2} p_j(u) X^j$, then $p_j(u), j=1,\ldots,q+1$ are linear polynomials with integer coefficients of $u$ and $p_0 = p_{q+2} =1$. Let $A_0(u), \ldots,A_{q+1}(u) \in \mathbb{Z}$ and define $A_{n}(u), n\ge q+2$ by the linear recursion 
\begin{equation} \label{e:ABC_rek}
A_{n}(u) = - p_{q+1}(u)A_{n-1}(u) -\ldots - p_0(u) A_{n-q-2}(u).
\end{equation}
For simplicity we denote the initial terms by $A_0,\ldots,A_{q+1}$, and to exclude the zero sequence we assume that at least one of them is not zero. Set $L= \max\{|A_1|,\ldots,|A_{q+1}|\}$.
 
Plainly $A_n(u)\in \mathbb{Z}[u]$ for all $n\ge 0$. The polynomial $P_u(X)$ is called the {\it characteristic polynomial} of the sequence $(A_n(u))$. As $p_0(u)=1$ this sequence can be extended to negative indexes too such that all members of the two sided sequence belong to $\mathbb{Z}[u]$. Hence the search for the integers, which appear multiple times in the two-parametric sequence $(A_n(u))$ is meaningful and non-trivial. In the language of Diophantine equations this means to find solutions of the equation
\begin{equation} \label{e:ABC}
|A_n(u)| = |A_m(u)|
\end{equation}
in $n,m,u \in \mathbb{Z}, n\not= m $. Our aim is to prove the following.

\begin{theorem} \label{t:ABC_gen}
Let $P_u(X)=XQ(X)(X-u)+1$, where $Q(X)\in \mathbb{Z}[X]$ is monic, separabel, has constant term $\pm 1$, and degree $q\ge 1$. Denote $(A_n(u))$ the linear recursive sequence having characteristic polynomial $P_u(X)$ and initial terms $A_0,\ldots,A_{q+1}\in \mathbb{Z}$ such that at least one of them is non-zero. 
Then there exist effectively computable constants $n_0,u_0$ depending only on $Q(X)$ and $L=\max\{|A_0|,\ldots, |A_{q+1}|\}$ such that if \eqref{e:ABC} holds for some integers $u,n,m$ with $u>u_0$ and $n\not=m$ then
$$
|n|,|m| <n_0.
$$
\end{theorem}

After submitting the manuscript we we made acquainted with the paper of Noubissie and Tichy \cite{NT}. They investigate common values of linear recursive sequences associated to generalized ABC-polynomials,  $\pi_U(X) = X (X-P_1(U))\cdots (X - P_{q-1}(U))+1$, where $P_1(U),\ldots, P_{q-1}(U) \in \mathbb{Z}[U]$ are pairwise distinct and such that $\deg P_1\le \ldots \le \deg P_{q-3}< \deg P_{q-2}< \deg P_{q-1}$. They proved that these sequences have at most two sufficiently large common values. Their result is not effective. Unfortunately it say nothing on common values of linear recursive sequences associated to our version of ABC-polynomials, but with $|Q(0)|>1$.

The strategy of our proof is similar that were used for the investigation of parametric families of Thue equations. In Section \ref{s:ABCpoly} we carefully investigate the roots of $P_u(X)$. It turns out that if $u$ is large enough then they behave regularly. It admits two real roots, one $\alpha_{q+1}(u)$ is as large as $u$, and the other $\alpha_0(u)$ lies near to zero, moreover $Q(0)\alpha_0(u) \alpha_{q+1}(u)$ is asymptotically one. 

Section \ref{s:ABCseq} is devoted to the ABC sequence with special attention to the weights of $\alpha_0(u)$ and $\alpha_{q+1}(u)$ appearing in the Binet formula of $(A_n(u))$. In both sections the machinery of Puiseux expansions of algebraic functions play important role. In Section \ref{s:transcendental} we collected necessary tools from transcendental number theory. 

In the last section we prove Theorem \ref{t:ABC_gen}. First we deal with the easy cases, i.e. if $n$ and $m$ have the same signs. The hard case is, when $nm<0$. Then we derive a lower bound for $\max\{|n|,|m|\}$ in term of $u$, which have the shape $g u^f$ with positive constants $g,f$. Having it we derive the desired upper bound for $\max\{|n|,|m|\}$ by using Matveev's theorem.

The symbols $c_x, n_x, u_x$ denote constants, which depend only on $Q(X)$ and $L$, the $u_x$ are such that the corresponding predicate holds only for $u\le u_x$, the $n_x$ play similar role with respect to $\max\{|n|,|m|\}$. 

\section{Tools from algebraic functions} 

To prove our main theorem in the generality as stated we have to invoke facts about algebraic functions. They are collected in this section. We present them for polynomials defined over $\mathbb{Q}$, because this is relevant in our application, although most results hold much more general. 

A function $S(z)\; : \; \mathbb{C} \mapsto \mathbb{C}$ is called {\it algebraic} if there is $F(X,Y)\in \mathbb{Q}[X,Y]$ such that $F(S(z),z)=0$. The algebraic functions form a field, see e.g. Jung \cite{Jung} (Viertes Kapitel, Satz I). The algebraic functions admit formal {\it Puiseux expansions}. More precisely let
$$
F(X,Y) = X^d+p_1(Y)X^{d-1} + \ldots + p_d(Y),
$$
where $p_1(Y),\ldots,p_d(Y)\in \mathbb{Q}[Y]$. The following assertion is the simplified version of Theorem 2.1.1 of \cite{EG}.

\begin{Pro} \label{p:Puiseux}
There are positive integers $r\le d$ and $e_1,\ldots,e_r$ with $e_1+\ldots+e_r = d$ and formal Puiseux series
$$
x_i(y) = \sum_{j=n_i}^{\infty} a_{ij}y^{-j/e_i},\; a_{in_i}\not=0,\; i=1,\ldots,r
$$
with $a_{ij}\in \mathbb{C}$ that satisfy $F(x_i(y),y)=0$. 
\end{Pro}

We need not only the formal Puiseux expansion of algebraic functions, but also their convergence if $y$ tends to infinity. Very special case of Lemma 1 of \cite{WMSchmidt} is that the Puiseux series of any algebraic functions have positive radius of convergence.

\section{Properties of the ABC polynomials}  \label{s:ABCpoly}

Although in Theorem \ref{t:ABC_gen} $Q(X)$ is assumed to have constant term $\pm1$, many preliminary statements hold without this assumption. In this section we suppose that $Q(X)\in \mathbb{Z}[X]$ is monic and separable, $Q(0)\not=0$ and that its degree is $q$. We set $P_u(X) = XQ(X)(X-u) + 1$, where $u$ denotes parameter, which is usually a rational integer, but sometimes real. In the later case this will be indicated explicitly. Then, plainly, $P_u(X) \in \mathbb{R}[X]$  ($P_u(X) \in \mathbb{Z}[X]$ if $u\in \mathbb{Z}$) is monic and of degree $q+2$. Let $\beta_0=0, \beta_{q+1} = u$ and denote $\beta_j,\; j=1,\ldots,q$ the zeroes of $Q(X)$. Then $\beta_j,\; j=0,\ldots,q+1$ are the zeroes of $XQ(X)(X-u)$ and they are pairwise distinct provided $u\not=\beta_j, j=1,\ldots,q$, which happens certainly if $|u|>\max \{|\beta_j|\;:\; j=1,\ldots,q\}$. Our first lemma generalizes Lemma 1 of \cite{ABC}. The proof is, up to notations, the same.

\begin{Le} \label{l:roots}
If $u\in \mathbb{R}$ and $u\ge u_1$ then $P_u(X)$ has the following properties:
\begin{enumerate}
\item[(i)] The degree of $P_u(X)$ is $q+2$, and its zeroes $\alpha_j(u), \; 0\le j\le q+1$ are pairwise distinct.
\item[(ii)] For a suitable arrangement of these roots 
$$
\lim _{u\to \infty} u (\alpha_i(u) - \beta_i) = b_i, \; i=0,\ldots,q
$$
hold with 
$$
b_0 = Q(0)^{-1}
$$
and
$$
b_i = \frac{1}{\beta_i Q'(\beta_i)}\not= 0, \; i=1,\ldots,q.
$$
\item[(iii)] If, with the arrangement of (ii), $\beta_i$ is real, then $\alpha_i(u)$ is real too. In particular the numbers $\alpha_0(u)$ and $\alpha_{q+1}(u)$ are real.  
\item[(iv)] If, additionally, $u\in \mathbb{Z}$ then $P_u(X)$ is irreducible in $\mathbb{Q}[X]$.
\item[(v)] If $u\in \mathbb{Z}$ then $\alpha_0(u)$ and $\alpha_{q+1}(u)$ are multiplicatively independent.
\end{enumerate}
\end{Le}

\begin{proof}
We do not separate the proofs of (i) and (ii). The degree of $P_u(X)$ is $q+2$, thus it has $q+2$ complex roots. We only have to prove that they are distinct. As $u\to \infty$ 
$$
-\frac{1}{u} P_u(X) \to  X Q(X)
$$
pointwise. Because of the continuity of the roots as functions of the parameter u, each of the $\beta_i,\; 0\le i\le q$, is the limit of a root $\alpha_i(u)$ of $P_u(X)$. Also, because of $P_u(\alpha_i(u)) = 0$, 
$$
(\alpha_i(u)-\beta_i) (\alpha_i(u)-u) \prod_{j=0 \atop j\not=i}^q (\alpha_i(u) - \beta_j)  +1 =0.
$$
Notice that 
$$
\lim_{u\to \infty} \prod_{j=0 \atop j\not=i}^q (\alpha_i(u) - \beta_j) = \prod_{j=0 \atop j\not=i}^q (\beta_i - \beta_j)^{} \not= 0.
$$
Hence, if $i>0$ then  
$$
\lim_{u\to \infty} u (\alpha_i(u) - \beta_i) = \prod_{j=0 \atop j\not=i}^q (\beta_i - \beta_j)^{-1} \not= 0,
$$
and for $i=0$ 
$$
\lim_{u\to \infty} u \alpha_0(u)  = \prod_{j=1 }^q ( - \beta_j)^{-1} = Q(0)^{-1} \not= 0,
$$
and all limits are non-zero constants. In particular, $\alpha_i(u) = \beta_i + O(u^{-1}),\; i=0,\ldots,q$. Since $\sum_{j=0}^{q+1} \alpha_j(u) = \sum_{j=0}^{q+1} \beta_j$, it follows that $\alpha_{q+1}(u) = \beta_{q+1} + O(u^{-1}) = u + O(u^{-1})$, which proves (ii).

To prove (iii) assume that $\beta_i$ is real, but $\alpha_i(u)$ is not real for some $0\le i\le q+1$ and for all large enough $u$. Then $b_i$ is real too. Further $\overline{\alpha_i(u)}$, the complex conjugate of $\alpha_i(u)$ is different from $\alpha_i(u)$, and is a zero of $P_u(X)$, i.e. $\overline{\alpha_i(u)}= \alpha_j(u)$ for some $j\not=i$. On the other hand  
$$
u(\alpha_j(u) -\beta_i) = u(\overline{\alpha_i(u)} -\beta_i) = \overline{u(\alpha_i(u)-\beta_i)} \to \overline{b_i}=b_i,
$$
but this contradicts (ii). 

Although the proof of (iv) repeats essentially the argumentation of Lemma 1 \cite{ABC} we present it for the convenience of the readers, and because our statement is more general. The assumption $u\in \mathbb{Z}$ yields $P_u(X)\in \mathbb{Z}[X]$. Assume that $P_u(X)$ is reducible in $\mathbb{Q}[X]$, i.e. $P_u(X)= g_1(X) g_2(X)$ with $g_1(X), g_2(X) \in \mathbb{Z}[X]$ both of degree at least one. We may assume that $\alpha_{q+1}(u)$ is a zero of $g_2(X)$. There is a $0\le j\le q$ such that $\alpha_j(u)$ is a zero of $g_1(X)$. Then $u -\alpha_j(u)$ is a unit because $P_u(u)=1=\prod_{i=0}^{q+1} (u-\alpha_i(u))$. Hence there is an algebraic conjugate of $u-\alpha_j(u)$, which absolute value is smaller than one. On the other hand the algebraic conjugates of $\alpha_j(u)$ are zeroes of $g_1(X)$, in particular different from $\alpha_{q+1}(u)$, thus by (ii) bounded. Hence for large enough $u$ all algebraic conjugates of $u-\alpha_j(u)$ have absolute value larger than one, which is a contradiction.

We come back to the proof of $(v)$ at the end of this section.
\end{proof}

\begin{Cor} \label{c:F_becslés}
Set $F=2\max\{|\beta_1|^{\pm1},\ldots,|\beta_q|^{\pm1}\}$. There exists a constant $u_2$ such that $\max\{|\alpha_1(u)|^{\pm1},\ldots,|\alpha_q(u)|^{\pm1}|\}< F$ and 
$$
|\alpha_0(u)|, \alpha_{q+1}(u)^{-1} < F <|\alpha_0(u)|^{-1}/2, \alpha_{q+1}(u)/2 
$$ 
hold for all $u\ge u_2$.
\end{Cor}

\begin{proof}
Both assertions follow from Lemma \ref{l:roots} (ii).
\end{proof}

\begin{Le} \label{l:szorzat}
Let $X, X_1,\ldots,X_q$ be indeterminates and set $Q(X) = \prod_{j=1 }^q (X-X_j)$. Denote $Q'(X)$ the derivative of $Q(X)$ with respect to $X$. Then
\begin{equation} \label{e:szorzat}
\sum_{j=1}^q \frac{1}{X_j^2 Q'(X_j)} = (-1)^{q-1}\frac{\sum_{i=1}^q \prod_{j=1  \atop j\not= i}^q X_j}{\prod_{j=1 }^q X_j^2} = \frac{Q'(0)}{Q(0)^2}.
\end{equation}

\end{Le}

\begin{proof}
Notice that 
$$
Q'(X) = \sum_{i=1}^q \prod_{j=1 \atop j\not= i}^q (X-X_j).
$$
This implies on one hand 
$$
Q'(0) =  (-1)^{q-1}\sum_{i=1}^q \prod_{j=1 \atop j\not= i}^q X_j
$$
and on the other hand
$$
Q'(X_j) = \prod_{i=1 \atop i\not= j}^q (X_j-X_i)
$$
for all $j=1,\ldots,q$. The first identity together with the Vieta relation
$$
Q(0) = (-1)^q \prod_{j=1 }^q X_j
$$
yields that the second equality in \eqref{e:szorzat} is correct. 

To prove the first equality we apply the usual addition rule and obtain
$$
\sum_{j=1}^q \frac{1}{X_j^2 Q'(X_j)} = \frac{N}{D},
$$
where 
$$
N = \sum_{i=1}^q \prod_{j=1 \atop j\not= i}^q \left(X_j^2 Q'(X_j)\right), \quad D= \left(\prod_{j=1}^q X_j^2  \right) \left(\prod_{j=1}^q Q'(X_j)\right).
$$
If $1\le i < j\le q$ then the multiplicity of $X_i-X_j$ is one in $Q'(X_i)$ and in $Q'(X_j)$, but zero in all other $Q'(X_h)$. Thus we have
$$
N = \sum_{h=1 \atop h\not= i,j}^q \prod_{k=1 \atop k\not=h}^q \left(X_k^2 Q'(X_k)\right) + (X_i-X_j) \left(X_i^2\prod_{h=1 \atop h\not= i,j}^q (X_i-X_h) - X_j^2\prod_{h=1 \atop h\not= i,j}^q (X_j-X_h)\right)V,
$$
because $X_i-X_j$ appears as $X_j-X_i = -(X_i-X_j)$ in $Q'(X_j)$. In the last expression 
$$
V= \prod_{h=1 \atop h\not=i,j}^q (X_h^2 Q'(X_h)).
$$
With the notation 
$$
T(X) = X^2 \prod_{h=1 \atop h\not= i,j}^q (X-X_h)
$$ 
we see that the expression in the large brackets is $T(X_i)-T(X_j)$, which is plainly divisible by $(X_i-X_j)$. For $1\le i_1<j_1\le q, 1\le i_2<j_2\le q$ and $(i_1,j_1)\not=(i_2,j_2)$ the polynomials $X_{i_1}-X_{j_1}$ and $X_{i_2}-X_{j_2}$ are coprime in $\mathbb{Z}[X_1,\ldots,X_q]$, thus the polynomial 
$$
U=\prod_{1\le i<j\le q} (X_i-X_j)^2
$$
divides the second summand of $N$, but it divides the first summand too, because in the product $\prod_{k=1 \atop k\not=h}^q (X_k^2 Q'(X_k)), h\not= i,j$ both $Q'(X_i)$ and $Q'(X_j)$ appear, thus the multiplicity of $X_i-X_j$ is two. Hence the  polynomial $U$ divides $N$, but it is, up to sign, equal to $\prod_{j=1}^q Q'(X_j)$, which is the second factor in $D$, thus it is a common factor of $N$ and $D$.

It remains to compute the quotient $N/U\in \mathbb{Z}[X_1,\ldots,X_q]$. Both polynomials $N$ and $U$ are symmetric in their variables, thus $N/U$ is a symmetric polynomial too. The degree of $X_i$ in $N$ is $2q-1$ because for $h\not=i$ the summand is $\prod_{j=1 \atop j\not= h}^q X_j^2Q'(X_j)$, i.e. $X_i$ appears $q+1$-times in $X_i^2 Q'(X_i)$ and one time in all other factors $X_j^2Q'(X_j)$, which number is $q-2$. On the other hand $X_i$ appears in exactly $2(q-1)$ factors of $U$, i.e. its degree in $U$ is $2q-2$. Hence $N/U$ is a symmetric polynomial, which is linear in all variables. 

Next we compare the maximum of the total degrees of the monomials in $N$ and $U$. The monomials in $U$ have the same total degree $q(q-1)$ because it is the product of $q(q-1)$ linear polynomials. The monomials in $X_j^2 Q'(X_j)$ have total degree $q-1+2 = q+1$. The product of $q-1$ such polynomials has total degree $(q-1)(q+1) = q^2-1$, which is an upper bound for the total degree of $N$. On the other hand the monomial $(X_2\cdots X_q)^{q+1}$ with total degree $q^2-1$ appears only in the first summand of $N$, hence it is a monomial in $N$. 

The consequence of the last two paragraphs is that the total degree of $N/U$ is $q^2-1 - (q^2-q) = q-1$, thus every monomial in $N/U$ is the product of $q-1$ indeterminates. On the other hand $N/U$ is linear in all indeterminates, i.e. 
$$
\frac{N}{U} = \sum_{i=1}^q \prod_{j=1 \atop j\not=i}^q X_j.
$$

To finish the proof we have to show that the sign on the right hand side of \eqref{e:szorzat} is correct. Multiplying \eqref{e:szorzat} by $\prod_{j=1}^q X_j^2$ the left hand side becomes
$$
\sum_{j=1}^q \frac{\prod_{i=1 \atop i\not=j}^q X_i^2}{Q'(X_j)}.
$$
Substituting $X_q=0$ only the term with $j=q$ remains because in all other summands $X_q$ appear as a factor. Further, as $Q'(X_q) = Q'(0) = (-1)^{q-1}\prod_{j=1}^{q-1} X_j$, we obtain that the left hand side is $ (-1)^{q-1}\prod_{j=1}^{q-1} X_j$. Substituting $X_q=0$ in the expression staying on the right hand side we obtain the same result, i.e. the sign is correct.
\end{proof}

\begin{Cor} \label{c:alfaszorzat}
We have 
$$
\alpha_1(u)\cdots \alpha_q(u) =  (-1)^q Q(0) - \frac{Q'(0)}{Q(0)u} + o\left(\frac{1}{u}\right).
$$
\end{Cor}

\begin{proof}
For any $\varepsilon>0$ there exists by Lemma \ref{l:roots} (ii) an $u_0(\varepsilon)$ such that for all $u>u_0(\varepsilon)$ there are $\delta_j$ with $|\delta_j|<\varepsilon $ such that
$$
\alpha_j(u) = \beta_j + \frac{b_j}{u} + \frac{\delta_j}{u}, \quad j=0,\ldots,q,
$$
with the $b_j$ given in the lemma. After multiplying for $j=1,\ldots,q$ the left and right hand sides of these equations we obtain
$$
\prod_{j=1}^q \alpha_j(u) = \prod_{j=1}^q \beta_j + \frac{1}{u} \left(\sum_{j=1}^q b_j\prod_{h=1\atop h\not=j}^q \beta_h\right) + \frac{\delta}{u},
$$
where $|\delta| \le \varepsilon H$, and $H$ is bounded by a quantity depending only on $q$ and $\beta_1, \ldots,\beta_q$. Using the expressions for the $b_i$ and that 
$$
(-1)^q\prod_{h=1}^q \beta_h = Q(0)
$$
we can rewrite the expression in the large brackets as
$$
(-1)^q Q(0) \sum_{j=1}^q \frac{1}{\beta_h^2 Q'(\beta_h)}.
$$
The sum is, by Lemma \ref{l:szorzat}, equal to $(-1)^{q-1} Q'(0)/Q^2(0)$, hence the last expression simplifies to $(-1)^q Q(0) (-1)^{q-1}  Q'(0)/Q^2(0) = -Q'(0)/Q(0)$. 
\end{proof}

\begin{Cor} \label{c:alfa3szorzat}
There exist integers $J,r>0$,  and $\pi_j,\mu_j\in \mathbb{C}, j\ge J, \pi_J\not=0$  such that
\begin{eqnarray*}
\frac1{Q(0)\alpha_0(u)\alpha_{q+1}(u)} &=& 1 + \frac{\pi_J}{u^{J/r}} + \sum_{j=J+1}^{\infty} \frac{\pi_j}{u^{j/r}}\\
Q(0)\alpha_0(u)\alpha_{q+1}(u) &=& 1 - \frac{\pi_J}{ u^{J/r}} + \sum_{j=J+1}^{\infty} \frac{\mu_j}{u^{j/r}}.
\end{eqnarray*}
Both series are convergent if $u>u_3$.
\end{Cor}

\begin{proof}
The identity $\alpha_0(u)\cdots \alpha_{q+1}(u) = (-1)^{q+2}= (-1)^q$ and Corollary \ref{c:alfaszorzat} yield
\begin{equation} \label{e:alpah_0alpha_q+1}
\frac{1}{\alpha_0(u)\alpha_{q+1}(u)} = \frac{\alpha_1(u) \cdots \alpha_q(u)}{(-1)^q } = Q(0) + \frac{(-1)^{q+1} Q'(0)}{Q(0)u} + o\left(\frac{1}{u} \right).
\end{equation}

The functions $\alpha_i(u),\; i=0,\ldots,q+1$ are algebraic, hence $\alpha_0(u)\alpha_{q+1}(u)$ and $1/\alpha_0(u)\alpha_{q+1}(u)$ are algebraic too, which are by Lemma \ref{l:roots} (iii) real for $u>u_0$. Hence, by Proposition \ref{p:Puiseux} there are $N\in \mathbb{Z}$, $\pi'_j\in \mathbb{R}, \; j\ge N, \pi'_J\not=0$ and  $1\le r\in \mathbb{Z}$ such that
$$
\frac1{\alpha_0(u)\alpha_{q+1}(u)} = \sum_{j=N}^{\infty} \pi'_j u^{-j/r}.
$$
Comparing this with \eqref{e:alpah_0alpha_q+1} we see that $N=0, \pi'_0=Q(0)$ and $\pi'_j=0$ for $j=1,\ldots,r-1$. Let $J$ be the smallest positive index with $\pi'_J\not=0$, then 
$$
\frac1{\alpha_0(u)\alpha_{q+1}(u)} = Q(0) + \frac{\pi'_J}{u^{J/r}} + \sum_{j=J+1}^{\infty} \frac{\pi'_j}{u^{j/r}}.
$$
After division by $Q(0)\not=0$ we get the first assertion with $\pi_j=\pi'_j/Q(0), \; j=0,1,\ldots$. If $Q'(0)\not=0$ then $J=r$ and $\pi_J=(-1)^{q+1}Q'(0)/Q(0)^2$. 

The algebraic function $Q(0)\alpha_0(u)\alpha_{q+1}(u)$ admits a Puiseux expansion too. The main term of its expansion is $1$. Assume that
$$
Q(0)\alpha_0(u)\alpha_{q+1}(u) = 1 + \frac{\mu_K}{u^{K/r}} + \sum_{j=K+1}^{\infty} \frac{\mu_j}{u^{j/r}} ,
$$
where $K$ denotes the smallest positive index with $\mu_K\not=0$. Multiplying the last two equations we get
\begin{eqnarray*}
1 &=& \left(1 + \frac{\pi_J}{u^{J/r}} + \ldots\right) \left( 1 + \frac{\mu_K}{u^{K/r}} + \ldots\right)\\
&=& 1 + \frac{\mu_K}{u^{K/r}}+  \frac{\pi_J}{u^{J/r}}\ldots,
\end{eqnarray*}
which is only possible, when $J=K$ and $\mu_J +\pi_J = 0$, i.e. if $\mu_J = -\pi_J$. The existence of the series is proved.

Their convergence is justified by Lemma 1 of Schmidt \cite{WMSchmidt}.
\end{proof}

The next corollary is a simplified, but more applicable version of Corollary \ref{c:alfa3szorzat}.

\begin{Cor} \label{c:alfa2szorzat}
There exist constants $\kappa_1\ge 1, \pi>0, u_4 $ depending only on $Q(X)$ such that $Q(0)\alpha_0(u)>0$ and either
$$
1 + \frac{\pi}{2u^{\kappa_1}} < \frac{1}{|Q(0)\alpha_0(u)| \alpha_{q+1}(u)} < 1 + \frac{2\pi}{u^{\kappa_1}}
$$ 
and 
$$
1 - \frac{2\pi}{u^{\kappa_1}} < |Q(0)\alpha_0(u)| \alpha_{q+1}(u)< 1 - \frac{\pi}{2u^{\kappa_1}} 
$$ 
or
$$
1 - \frac{2\pi}{u^{\kappa_1}} < \frac{1}{|Q(0)\alpha_0(u)| \alpha_{q+1}(u)} < 1 - \frac{\pi}{2u^{\kappa_1}}
$$ 
and 
$$
1 + \frac{\pi}{2u^{\kappa_1}} < |Q(0)\alpha_0(u)| \alpha_{q+1}(u)< 1 + \frac{2\pi}{u^{\kappa_1}} 
$$ 
hold for all $u>u_4$.
\end{Cor}

\begin{proof}
An obvious consequence of Corollary \ref{c:alfa3szorzat} is that $Q(0)\alpha_0(u)\alpha_{q+1}(u)>0$. As $\alpha_{q+1}(u)>0$ we get $Q(0)\alpha_0(u)>0$.

By the second assertion of Corollary \ref{c:alfa3szorzat} the Puiseux series of $\alpha_0(u)\alpha_{q+1}(u)$ and of $1/\alpha_0(u)\alpha_{q+1}(u)$ are convergent, whenever $u>u_0$. In particular 
$$
\left|\sum_{j=J+1}^{\infty} \frac{\pi_j}{u^{j/r}}\right| < \frac{|\pi_J|}{2u^{J/r}} \quad \mbox{and} \quad \left|\sum_{j=J+1}^{\infty} \frac{\mu_j}{u^{j/r}}\right| < \frac{|\mu_J|}{2u^{J/r}} = \frac{|\pi_J|}{2|Q(0)|u^{J/r}}
$$
hold for all $u>u_0$.

If $\pi_J>0$ then $1/|Q(0)\alpha_0(u)|\alpha_{q+1}(u)> 1$ and we obtain the first alternative with $\pi = \pi_J$ and $\kappa_1=J/r$. The opposite case yields the second alternative, with $\pi = -\pi_J$ and $\kappa_1=J/r$.
\end{proof}

\medskip
Now we are enough equipped to prove Lemma \ref{l:roots} $(v)$.

\medskip
{\it Proof of Lemma \ref{l:roots} $(v)$.} 
Assume that $u\ge u_1$ and $\alpha_0(u)$, $\alpha_{q+1}(u)$ are multiplicatively dependent, i.e. there exist integers $a,b$ at least one of them non-zero, such that. $\alpha_0^a(u) \alpha_{q+1}^b(u) = 1$. If $b<0$ then $\alpha_0^{-a}(u) \alpha_{q+1}^{-b}(u) = 1$, which is an analogous equation, but with $-b>0$. Hence we may assume $b\ge 0$ in the sequel. By Lemma \ref{l:roots} $(iv)$ the polynomial $P_u(X)$ is irreducible. If $b=0$ then $\alpha_0(u)$ is a root of unity, but this is impossible by the irreducibility of $P_u(X)$. Similarly $a\not=0$. Hence $b>0$ and $a\not=0$.

The case $a<0, b>0$ is impossible, because $|\alpha_0(u)|<1$ and $|\alpha_{q+1}(u)|>1$, thus $|\alpha_0^a(u) \alpha_{q+1}^b(u)| > 1$, contradicting the assumption.

If $a=b$ then $(\alpha_0(u) \alpha_{q+1}(u))^a = 1$, which yields $|\alpha_0(u) \alpha_{q+1}(u)| = 1$, but this contradicts Corollary \ref{c:alfa2szorzat}. 

Denote $\Sigma$ the Galois group of the normal closure of $\mathbb{Q}(\alpha_0(u))$. If there is an element $\sigma \in \Sigma$, which swaps $\alpha_{q+1}(u)$ and $\alpha_{0}(u)$, then 
$$
1=\sigma(\alpha_0^a(u) \alpha_{q+1}^b(u)) = \sigma(\alpha_0(u))^a \sigma(\alpha_{q+1}(u))^b = \alpha_{q+1}^a(u) \alpha_0^b(u).
$$
This together with the original equation yields
$$
(\alpha_0(u) \alpha_{q+1}(u))^{a+b} =1,
$$
which was excluded earlier. 

It remains to deal with the case when $a,b>0$ and no elements of $\Sigma$ swaps $\alpha_{q+1}(u)$ and $\alpha_{0}(u)$. We distinguish two cases:
$0<a<b$ and $0<b<a$.

We assume first $0<a<b$. As the size of $\Sigma$ is at least $q+2\ge3$ there is an element $\sigma \in \Sigma$, for which $\sigma(\alpha_{q+1}(u)) = \alpha_{0}(u)$ and $\sigma(\alpha_{0}(u)) = \alpha_{j}(u), j\not=0,q+1$. Now we have
$$
1=\sigma(\alpha_0^a(u) \alpha_{q+1}^b(u)) = \alpha_j^a(u) \alpha_{0}^b(u) = (\alpha_j(u) \alpha_{0}(u))^a \alpha_{}^{b-a}(u),
$$
but this is impossible because $|\alpha_0(u) \alpha_{j}(u)|<1$ and $|\alpha_{0}(u)|<1$ whenever $u$ is large enough.

In the second case, i.e. if $0<b<a$ we choose a $\sigma \in \Sigma$, for which $\sigma(\alpha_{0}(u)) = \alpha_{q+1}(u)$, and $\sigma(\alpha_{q+1}(u)) = \alpha_{j}(u), j\not=0,q+1$. Now we have
$$
1=\sigma(\alpha_0^a(u) \alpha_{q+1}^b(u)) = \alpha_{q+1}^a(u) \alpha_{j}^b(u) = \alpha_{q+1}^{a-b}(u)(\alpha_0(u) \alpha_{j}(u))^b,
$$
but this is impossible because $|\alpha_j(u) \alpha_{q+1}(u)|>1$ and $|\alpha_{q+1}(u)|>1$ whenever $u$ is large enough.$\Box$

\section{Properties of the ABC sequences}  \label{s:ABCseq}

Now we turn to the study of $(A_n(u))$. We assume throughout this section that at least one of $A_0,\ldots,A_{q+1}$ is non-zero. In the Introduction, more precisely by \eqref{e:ABC_rek} the sequence $(A_n(u))$ is defined for non-negative indices, but we pointed out that it can be extended uniquely for negatives too, such that all members of the two-directional sequence are polynomials with integer coefficients. This is a consequence of the next lemma, which shows that $(A_n(u))$ is a linear recursive sequence with integer coefficients for negative indices too.

\begin{Le} \label{l:negative_rekurzio}
If the characteristic polynomial of the sequence $(A_n(u)), \; n\ge 0$ is $P_u(X)$ then it is a linear recursive sequence for $n<0$ too with the characteristic polynomial $X^{q+2}P_u(1/X)$.
\end{Le}

\begin{proof}
Set $p_0(u)=p_{q+2}(u) = 1$ and add $p_0(u) A_{n-q-2}(u) - A_n(u) = A_{n-q-2}(u) - p_{q+2}(u)A_n(u) $ to both sides of \eqref{e:ABC_rek}, then, after rearrangement of the RHS, we get
$$
A_{n-q-2}(u) = -p_1(u)A_{n-q-1}(u) - \ldots - p_{q+1}(u) A_{n-1}(u) - p_{q+2}(u) A_n(u).
$$
The substitution $n-q-2=-m$ yields
$$
A_{-m}(u) = -p_1(u)A_{-(m-1)}(u) - \ldots - p_{q+1}(u) A_{-(m-q-1)}(u) - p_{q+2}(u) A_{-(m-q-2)}(u),
$$
and the lemma is proved.
\end{proof}

Although the members of $(A_n(u))$ satisfy a linear relation, they grow typically exponentially. The Binet formula express this unexpected phenomena, see e.g. Theorem C.1. in \cite{ShT}. Assume that $u\ge u_1$ then (i) of Lemma \ref{l:roots} holds, i.e. $\alpha_0(u),\ldots,\alpha_{q+1}(u)$ are pairwise different. By the just quoted theorem there exist uniquely determined functions $a_0(u),\ldots,a_{q+1}(u)$, which belong for any fixed $u$ to the algebraic number field $\mathbb{Q}(\alpha_0(u),\ldots,\alpha_{q+1}(u))$ such that 
\begin{equation} \label{e:Binet}
A_n(u) = a_0(u)\alpha_0(u)^n + \ldots + a_{q+1}(u) \alpha_{q+1}(u)^n
\end{equation}
holds for all $n\in \mathbb{Z}$. The next lemma describes properties of the weights in the Binet formula. Denote $({\bf b_1},\ldots,{\bf b_d} )$ the matrix having columns ${\bf b_1},\ldots,{\bf b_d} $. Further set 
\begin{eqnarray*}
{\bf A} &=&(A_0,\ldots,A_{q+1})^t,\\ 
{\bf v}_j&=&(\alpha_j(u)^0,\ldots,\alpha_j(u)^{q+1})^t, \; j=0,\ldots,q+1,\\
D(u) &=& \det ({\bf v}_0,\ldots,{\bf v}_{q+1}) \quad \mbox{and}\\
D_j(u) &=& \det ({\bf v}_0,\ldots,{\bf v}_{j-1}, {\bf A}, {\bf v}_{j+1},\ldots,{\bf v}_{q+1}), \; j=0,\ldots,q+1.
\end{eqnarray*}

\begin{Le} \label{l:a0:a_q+2}
The following assertions are true for all $u_1<u\in \mathbb{Z}$
\begin{itemize}
\item[(i)] $a_j(u) \not= 0$ for all $j=0,\dots,q+1$,
\item[(ii)] If for some $0\le j\le q+1$ the function $\alpha_j(u)$ is real, then $a_j(u)$ is real too. Especially $a_0(u)$ and $a_{q+1}(u)$ are real,
\item[(iii)] 
$$
\frac{a_{q+1}(u)}{a_{0}(u)} = \frac{D_{q+1}(u)}{D_0(u)}.
$$
\item[(iii)] Denote $\Sigma$ the Galois group of the normal closure of $\mathbb{Q}(\alpha_0(u))$. If $\sigma\in \Sigma$ swaps $\alpha_i(u)$ and $\alpha_j(u)$ for some $0\le i\le j\le q+1$ then $\sigma(D_j(u)) = \pm D_i(u)$.
\end{itemize}
\end{Le}

\begin{proof}
Substituting $n=0,1,\ldots,q+1$ into \eqref{e:Binet} we obtain the system of linear equations
$$
A_j(u) = a_0(u)\alpha_0(u)^j +\ldots + a_{q+1}(u) \alpha_{q+1}(u)^j, \; j=0\ldots,q+1
$$
for the unknowns $a_0(u),\ldots,a_{q+1}(u)$. The coefficients form a Vandermonde matrix, whose determinant is $D(u) = \prod_{0\le i<j\le q+1} (\alpha_i(u) - \alpha_j(u))$, which is non-zero by Lemma \ref{l:roots} (i). Hence the above system has unique solution. Cramer's rule tells us, that
\begin{equation} \label{e:aj}
a_j(u) = \frac{D_j(u)}{D(u)}.
\end{equation}

By Lemma \ref{l:roots} (iv) $P_u(X)$ is irreducible, hence the $\alpha(u)$'s are algebraic conjugates. Thus the interchange of $\alpha_j(u)$ and $\alpha_0(u)$ is an isomorphism between the number fields $\mathbb{Q}(\alpha_j(u))$ and $\mathbb{Q}(\alpha_0(u))$, which does not affect $\mathbb{Q}$. Hence, if $D_j(u)=0$ for some $0\le j \le q+1$, then $D_j(u)=0$ would hold for all $0\le j \le q+1$, i.e. $A_n(u)=0$ for all $n\ge q+1$. This together with Lemma \ref{l:negative_rekurzio} implies $A_n(u)=0$ for all $n\ge 0$, but this contradicts the assumption regarding the initial terms. Hence (i) is true.

To prove (ii) observe that we obtain the matrix of $D_j(u)$ from the matrix of $D(u)$ by replacing the real column vector ${\bf v}_j$ by an other real column vector $\bf A$. The complex conjugation implies the interchange of the column vectors corresponding conjugate complex roots of $P_u(X)$, but the total number of interchanges are the same in both determinants, hence this transformation does not affect $a_j(u)$, i.e. it is a real number, and (ii) is proved.

Assertion (iii) is a direct consequence of $\eqref{e:aj}$ and $D_0(u)\not= 0$.

Now we prove $(iv)$. By definition  
$$
D_j(u) = \det ({\bf v}_0,\ldots,{\bf v}_{j-1}, {\bf A}, {\bf v}_{j+1},\ldots,{\bf v}_{q+1}), \; j=0,\ldots,q+1.
$$
Let $\sigma$ denote arbitrary element of $\Sigma$. If $0\le k\le q+1$ and $\sigma(\alpha_{k}(u)) = \alpha_{l}(u)$ then 
$$
\sigma({\bf v}_k) = \sigma((1, \alpha_{k}(u),\dots, \alpha_{k}^{q+1}(u))) = (1, \alpha_{l}(u),\dots, \alpha_{l}^{q+1}(u)) = {\bf v}_l,
$$
hence
$$
\sigma(D_j(u)) = \det (\sigma({\bf v}_0),\ldots,\sigma({\bf v}_{j-1}), \sigma({\bf A}), \sigma({\bf v}_{j+1}),\ldots,\sigma({\bf v}_{q+1})).
$$
Notice that $\sigma({\bf A}) = {\bf A}$, because $\sigma$ is the identity mapping on $\mathbb{Q}$. 

After this preparation we prove the assertion. If $j=i$ then $\sigma$ is a permutation on the set $\{{\bf v}_0,\ldots,{\bf v}_{j-1}, {\bf v}_{j+1},\ldots,{\bf v}_{q+1}\}$, hence the column vectors of $\sigma(D_j(u))$ and $D_j(u)$ are up to permutation identical, and the claim is proved in this case.

If $j\not=i$, then the $i$-th column of the determinant $\sigma(D_j(u))$ changes to ${\bf v}_j$, the $j$-th column does not change, while the other columns are a permutation of the vectors  ${\bf v}_0,\ldots,{\bf v}_{j-1}, {\bf v}_{j+1},\ldots,{\bf v}_{q+1}$, hence $\sigma(D_j(u))$ is up to permutation of its columns identical to $D_i$, and the claim is proved in this case too.
\end{proof}

In the next step we prove lower and upper bounds for $|D(u)|, |D_0(u)|$ and $|D_{q+1}(u)|$.

\begin{Le} \label{l:Dbounds}
There exist easily computable constants $c_1,c_2,c_3,c_4$ depending only on the degree and the roots of $Q(X)$ such that
\begin{eqnarray}
c_1 \alpha_{q+1}(u)^{q+1} &\le& |D(u)| \le c_2 \alpha_{q+1}(u)^{q+1} \label{e:D} \\
\frac1{(c_3 c_4^{q-1}L^q \alpha_{q+1}(u)^{q^2-1})}  &\le& |D_j(u)| \le c_4 L \alpha_{q+1}(u)^{q+1}, \quad j=0,\ldots,q \label{e:D0} \\
\frac1{( (c_4 L)^q \alpha_{q+1}(u)^{q(q+1)})} &\le& |D_{q+1}(u)| \le c_3 L \label{e:Dq}
\end{eqnarray}
hold for all $u\in \mathbb{Z}$, provided $P_u(X)$ is irreducible.
\end{Le}

\begin{rem}
If $u \ge u_1$ then $P_u(X)$ is irreducible by (iv) of Lemma \ref{l:roots}.
\end{rem}

\begin{proof}
If $P_u(X)$ is irreducible, then $\alpha_0(u),\ldots,\alpha_{q+1}(u)$ are pairwise distinct complex numbers, and by Lemma \ref{l:a0:a_q+2} the determinants $D(u), D_j(u), \; j=0,\ldots,q+1$ are non-zero. 

We know that $D(u)$ is a Vandermonde determinant, more exactly 
$$
D(u) = \prod_{0\le i< j \le q+1} (\alpha_i(u)-\alpha_j(u)),
$$
hence 
$$
|D(u)| = \alpha_{q+1}(u)^{q+1} \prod_{j=0}^q \left|1-\frac{\alpha_j(u)}{\alpha_{q+1}(u)} \right| \prod_{0\le i< j \le q} |\alpha_i(u)-\alpha_j(u)|,
$$
and \eqref{e:D} follows by Lemma \ref{l:roots}.

To prove the other two assertions we need the inequality
\begin{equation} \label{e:szorzat1}
\left|\prod_{j=0}^{q+1}D_j(u)\right| \ge 1
\end{equation}
for all appropriate $u$. To prove this pick an $u\ge u_1$ and show that $\pm D_j(u),\; j=0,\ldots,q+1$ are all the algebraic conjugates of $D_0(u)$. If this is true, then $|\prod_{j=0}^{q+1}D_j(u)|$ is a rational number, but $D_0(u)$ is an algebraic integer, thus the above product, which is the absolute value of the norm of $D_0(u)$, is a rational integer. The $D_0(u)$ is an algebraic integer, because all entries of the matrix defined it are algebraic integers. On the other hand we showed in the last lemma that $D_0(u)\not= 0$, which yields the claim. 

It remains to show that $\pm D_j(u), j=0,\ldots,q+1$ are all the algebraic conjugates of $D_0(u)$. Indeed, let $\tau$ be an algebraic conjugation of $\alpha_0(u)$, which can be interpreted as a permutation of the set $\{\alpha_0(u),\ldots,\alpha_{q+1}(u)\}$. Each permutation is the composition of transpositions, which are mappings interchanging two elements. Assume that $\tau$ is the transposition $\tau(\alpha_i(u))= \alpha_j(u)$ and $\tau(\alpha_j(u))= \alpha_i(u)$ for some $0\le i \le j \le q+1$. Restricting $\tau$ to $\mathbb{Q}$ is the identity mapping, thus $\tau({\bf A})= {\bf A}$.

If $i=j$ then $\tau$ is the identity mapping, in particular $\tau(D_0)=D_0$. If $i>0$ then $\tau$ interchanges the vectors ${\bf v}_i$ and ${\bf v}_j$ in $D_0$, hence $\tau(D_0)=\pm D_0$. Finally, if $0=i<j$, then $\tau({\bf v}_j) = {\bf v}_0$, hence $\tau(D_0)=\det({\bf A},{\bf v}_1,\ldots,{\bf v}_{j-1},{\bf v}_0,{\bf v}_{j+1}, \ldots,{\bf v}_{q+1})$, but this is $\pm D_j$, and our claim is proved.

Let $0\le j \le q$. Then, by the Hadamard inequality, we have
$$
|D_j(u)| \le |{\bf A}| \cdot|{\bf v}_{q+1}|\prod_{k=0 \atop k\not=j}^{q}|{\bf v}_k| \le (q+1)L\alpha_{q+1}(u)^{q+1} \prod_{k=0 \atop k\not=j}^{q}|{\bf v}_k|,
$$
where $|{\bf v}|$ denotes the Euclidean length of the vector ${\bf v}$. Corollary  \ref{c:F_becslés} implies that $((q+1)F)^{q-1}$ is an upper bound for the third factor, which depends only on the roots and degree of $Q(X)$. This proves the upper bound in \eqref{e:D0} with 
$$
c_4 = ((q+1)F)^{q}.
$$

The upper bound for $|D_{q+1}(u)|$ can be proved analogously by observing that $|{\bf v}_0|$ has an absolute upper bound. The lower bounds follow now by combining the upper bounds and \eqref{e:szorzat1}.
\end{proof}

In the proof of our theorem we need not only the quotient $a_{q+1}(u)/a_0(u)$, but its limit when $u$ tends to infinity too. The functions $\alpha_0(u),\ldots,\alpha_{q+1}(u)$ are by Lemma \ref{l:roots} (i) pairwise different for all $u_1\le u\in \mathbb{R}$. They yield the functions $D(u), D_j(u) \; :\; \mathbb{R} \mapsto \mathbb{C}, \; j=0,\ldots,q+1$. As $D(u)\not=0$ whenever $u\ge u_1$, the functions $a_j(u) = \frac{D_j(u)}{D(u)} \; :\; \mathbb{R} \mapsto \mathbb{C}, \; j=0,\ldots,q+1$ are well defined too. Unfortunately the limit of the quotient 
$$
\frac{a_{q+1}(u)}{a_0(u)} = \frac{D_{q+1}(u)}{D_0(u)} 
$$
is not always well defined or is zero. Obvious instance of the later is 
because $D_0(u)$ tends with $\alpha_{q+1}(u)$ typically to infinity. Hence it is more convenient to study the limit of $D_0(u)/\alpha_{q+1}(u)^{q+1}$, which we denote by ${\mathcal D}_0'$. Lemma \ref{l:roots} (ii) yields
$$
{\mathcal D}_{0}' = \left| 
\begin{array}{ccccc}
A_0&1&\ldots&1&0\\
A_1 &\beta_1&\ldots&\beta_q&0\\
\ldots& \ldots& \ldots& \ldots& \ldots\\
A_q&\beta_1^{q}&\ldots&\beta_q^{q}&0\\
A_{q+1}&\beta_1^{q+1}&\ldots&\beta_q^{q+1}&1\\
\end{array}
\right| = \left| 
\begin{array}{cccc}
A_0&1&\ldots&1\\
A_1&\beta_1&\ldots&\beta_q\\
\ldots& \ldots& \ldots& \ldots\\
A_q&\beta_1^{q}&\ldots&\beta_q^{q}\\
\end{array}
\right|.
$$
Denote $d_0$ the determinant staying on the right hand side.

Set further
$$
{\mathcal D}_{q+1}' = \lim_{u \to \infty} D_{q+1} (u).
$$
Lemma \ref{l:roots} (ii) yields
$$
{\mathcal D}_{q+1}' = \left| 
\begin{array}{ccccc}
1&1&\ldots&1&A_0\\
0 &\beta_1&\ldots&\beta_q&A_1\\
\ldots& \ldots& \ldots& \ldots& \ldots\\
0 &\beta_1^{q+1}&\ldots&\beta_q^{q+1}&A_{q+1}\\
\end{array}
\right| = \beta_1\cdots \beta_q \left| 
\begin{array}{cccc}
1&\ldots&1&A_1\\
\beta_1&\ldots&\beta_q&A_2\\
\ldots& \ldots& \ldots& \ldots\\
\beta_1^{q}&\ldots&\beta_q^{q}&A_{q+1}\\
\end{array}
\right|.
$$
Setting 
$$
d_{q+1}= \left| 
\begin{array}{cccc}
1&\ldots&1&A_1\\
\beta_1&\ldots&\beta_q&A_2\\
\ldots& \ldots& \ldots& \ldots\\
\beta_1^{q}&\ldots&\beta_q^{q}&A_{q+1}\\
\end{array}
\right|
$$
we have ${\mathcal D}_{q+1}'=(-1)^qQ(0)d_{q+1}$.

Unfortunately it can happen that $d_0$ or $d_{q+1}$ or both are zero. For example if $A_0=\ldots = A_{q}=0, A_{q+1}=1$ then, obviously, $d_0=0$. Another example is when $\beta_1=1$ and $A_j=1,\; j=0,\ldots,q+1$, in which case $d_0=d_{q+1}=0$. I did several, attempts to characterize that irregular cases. Although I was not successful, it was possible to prove Theorem \ref{t:ABC_gen} under stronger assumptions, but with explicit constants. This can be the subject of a different paper, because the technical details would extend the size of these notes considerably. To overcome this difficulty we use that our determinants $D_0(u), D_{q+1}(u)$ are algebraic functions, which makes it possible to deploy the machinery of that theory.

\begin{Le} \label{l:Dq+1/D0v2}
There exist  effective constants $\kappa,\kappa_2\in \mathbb{Q}, \varphi,\mu, u_5\in \mathbb{R}, |\varphi|,\mu >0$ depending only on $Q,A_0,\ldots,A_{q+1}$ such that the function 
$$
\Phi(u)= u^{\kappa} \frac{D_{q+1}(u)}{D_{0}(u)}
$$
satisfy
$$
\lim_{u\to \infty} \Phi(u)= \varphi.
$$
Moreover, if $ \Phi(u)$ is not ultimately constant, then 
$$
\frac{\mu}{2u^{\kappa_2}}<| \Phi(u)-\varphi| < \frac{2\mu}{u^{\kappa_2}}
$$
holds for all $u\ge u_5, u\in \mathbb{Z}$.
\end{Le}

\begin{proof}
The functions $D_0(u), D_{q+1}(u)$ are the values of determinants whose entries are the integers $A_0,\ldots,A_{q+1}$ and the powers of $\alpha_0(u),\ldots,\alpha_{q+1}(u)$, which are the roots of $P_u(X)\in \mathbb{Z}[X]$, hence algebraic functions. Thus $D_{0}(u), D_{q+1}(u)$ are algebraic functions, consequently they admit Puiseux expansions. By Proposition \ref{p:Puiseux} there are integers $n_i, i=0,q+1$, positive integers $r_i, i=0,q+1$ and $a_{ij}\in \mathbb{C},\; i=0,q+1, \; j\ge n_i$ such that
$$
D_i(u) = \sum_{j=n_i}^{\infty} a_{ij} u^{-j/r_i}, \;a_{in_i}\not=0
, \; i=0,q+1.
$$
If $r_0\not=r_{q+1}$ then replacing $n_i$ by $n_ir/r_i, \; i=0,q+1$, where $r$ denotes the least common multiple of $r_0,r_{q+1}$, the transformed Puiseux series have the property $r_0=r_{q+1}$. In the rest of this proof we write $r$ instead of $r_0$ or $r_{q+1}$.

Lemma \ref{l:Dq+1/D0v2} ensures that 
\[
|D_0(u)| \alpha_{q+1}^{q^2-1} \quad \mbox{is larger, and} \quad \frac{|D_{0}(u)|} {\alpha_{q+1}^{q+1}} \quad \mbox{is less, than some positive constants for all } 
\]
large enough $u$, hence 
$$
-(q+1) \le \frac{n_0}{r} \le q^2-1.
$$
We used here that $\alpha_{q+1}(u)/u$ tends to one very rapidly. Similar argumentation leads to
$$
0 \le \frac{n_{q+1}}{r} \le q(q+1).
$$

With $D_0(u)$ and $D_{q+1}(u)$ their quotient is an algebraic function having formal Puiseux expansion, i.e. there exist integers $J,r$ and non-zero complex numbers $\pi_J$ such that
$$
\frac{D_{q+1}(u)}{D_{0}(u)} = \sum_{j=J}^{\infty} \pi_{j} u^{-j/r}.
$$
The inequalities for $n_i/r, \; i=0, q+1$ imply 
$$
-q^2+1 \le \frac{J}{r} \le (q+1)^2.
$$

The choice $\kappa=J/r, \varphi=\pi_J$ proves the first statement, because  
$$
\Phi(u) = u^{\kappa} \frac{D_{q+1}(u)}{D_u(u)} = \varphi +\sum_{j=1}^{\infty} \pi_{J+j} u^{-j/r}.
$$
If $\Phi(u)$ is not ultimately constant, then $\Phi(u)\not= \varphi$ for all $u\ge u_5$, but then there exists $j>0$ with $\pi_{J+j}\not=0$ too. Choosing $J_0$ the smallest such index, the assertion holds with $\mu= |\pi_{J+J_0}|, \kappa_2=J_0/r$.

\end{proof}

To prove our theorem we need estimations on the growth of $|A_n(u)| = |A_n({\bf A},u)|$ and $|A_{-n}(u)| = |A_{-n}({\bf A},u)|$.

\begin{Le} \label{l:novekedes}
Let $F=2\max\{|\beta_1|^{\pm1},\ldots,|\beta_q|^{\pm1}\}$, like in Corollary \ref{c:F_becslés}. There exists an effectively computable constant $n_2$, depending only on $Q(X)$ and $L$, such that if $n\ge n_2$ then
\begin{equation} \label{e:pozitiv}
||A_{n}(u)|-|a_{q+1}(u)\alpha_{q+1}(u)^{n}| |< c_8 L F^n
\end{equation} 
and 
\begin{equation} \label{e:negativ}
||A_{-n}(u)|-|a_0(u)\alpha_0(u)^{-n}| |< c_8 L F^n
\end{equation} 
hold for all $u\ge u_2$.
\end{Le}

\begin{proof}
By the Binet formula 
$$
A_{m}(u) = \sum_{j=0}^{q+1} a_j(u) \alpha_j(u)^m
$$
holds for every integer $m$. 

Suppose that $m=n$ with a non-negative integer $n$, then 
\begin{eqnarray*}
|A_{n}(u)-a_{q+1}(u)\alpha_{q+1}(u)^{n}| &=& \left|\sum_{j=0}^{q} a_j(u) \alpha_j(u)^n \right|\\
&\le & (q+1) c_4L \max\{|\alpha_j(u)|^n,\; j=0,\ldots,q\}/c_1 \\
&\le & (q+1) c_4L F^n/c_1.
\end{eqnarray*}
The first inequality holds because \eqref{e:D} and \eqref{e:D0}, while the second because $|\alpha_j(u)|<F, \; j=0,\ldots,q$ by Corollary \ref{c:F_becslés}.

If $A_{n}(u)$ and $a_{q+1}(u)\alpha_{q+1}(u)^{n}$ have the same signs, then \eqref{e:pozitiv} holds immediately with $c_8 = (q+1)c_4/c_1$. Otherwise   
$$
|a_{q+1}(u)\alpha_{q+1}(u)^{n}|< (q+1) c_4L F^n/c_1.
$$
Now, by \eqref{e:D} and \eqref{e:Dq} we have
\begin{eqnarray*}
|a_{q+1}(u)| &=& \frac{|D_{q+1}(u)|}{|D(u)|} \\
&>& \frac{1}{c_4L \alpha_{q+1}(u)^{q+1}} \frac{1}{c_2 \alpha_{q+1}(u)^{q+1}}\\
&=& \frac{1}{c_2 c_4L \alpha_{q+1}(u)^{2(q+1)}}.
\end{eqnarray*}
Combining the last two inequalities together we obtain
$$
2^n<\left(\frac{\alpha_{q+1}(u)}{F}\right)^n < \frac{(q+1)c_2(c_4L)^2}{c_1} |\alpha_{q+1}(u)^{2(q+1)}|,
$$
which yields an easily computable upper bound $n_2'$ for $n$. The first inequality holds because $\alpha_{q+1}(u)/F>2$ by Corollary \ref{c:F_becslés}.

If $m=-n$ with a non-negative integer $n$ then the dominating term in the Binet formula is $a_0(u)\alpha_0^{-n}$, and, after this observation, the proof is identical to the first case.
\end{proof}

\begin{Cor} \label{c:monoton}
If $u>u_2$ then $(|A_n(u)|)$ and $(|A_{-n}(u)|)$ are strictly monotone increasing whenever $n\ge n_3$. 
\end{Cor}

\begin{proof}
If $n>n_2 $ then by Lemma \ref{l:novekedes} we have 
$$
|A_{n+1}(u)| > |a_{q+1}(u)\alpha_{q+1}(u)^{n+1}| - c_8LF^{n+1}
$$
and 
$$
|A_{n}(u)| < |a_{q+1}(u)\alpha_{q+1}(u)^{n}| + c_8LF^{n}.
$$
If for such an $n$ the inequality $|A_{n+1}(u)| \le |A_{n}(u)|$ holds then
$$
|a_{q+1}(u)\alpha_{q+1}(u)^{n+1}| - c_8LF^{n+1} < |a_{q+1}(u)\alpha_{q+1}(u)^{n}| + c_8LF^{n}.
$$
Rearranging this inequality we obtain
$$
|a_{q+1}(u)\alpha_{q+1}(u)^{n}|( \alpha_{q+1}(u) -1) < c_8LF^{n} (F+1),
$$
but this yields $n\le n_3$ because $\alpha_{q+1}(u) > 2F$.

The proof that $(|A_{-n}(u)|)$ are strictly monotone increasing is the same.
\end{proof}

\section{Tools from transcendental number theory} \label{s:transcendental}

An important ingredient of the proof of the main theorem is A. Baker's theory on linear forms in logarithms of algebraic numbers. It has many variants, nowadays one of the sharpest is due to Matveev \cite{M}.


For an algebraic number $\eta$ with minimal polynomial
$$
f(X)=a_0(X-\eta^{(1)})\cdots(X-\eta^{(N)})\in {\mathbb Z}[X]
$$
with $a_0>0$ and with relative prime coefficients, write $h(\eta)$ for its {\it absolute logarithmic or Weil height} given by
$$
h(\eta)=\frac{1}{N}\left(\log a_0+\sum_{j=1}^N \max\{0,\log |\eta^{(j)}|\}\right).
$$
If $\eta$ is an algebraic integer then
\begin{equation} \label{e:heta}
h(\eta) \le \log \max\{|\eta^{(j)}|,\; j=1,\ldots,N\}.
\end{equation}

We provide some important and well known properties of the function $h$.


\begin{Le}\label{l:Weil height}
Let $\gamma,\eta$ be algebraic numbers of degree at most $d$ and $u\in \Q$. Then we have
\begin{enumerate}
  \item $h(\gamma \pm \eta) \le h(\gamma)+ h(\eta)+ \log 2$,
  \item $h(\gamma \eta^{\pm1}) \le h(\gamma)+ h(\eta)$,
  \item $h(\gamma^u) = |u| h(\gamma)$.
\end{enumerate}
\end{Le}
For the proof see e.g. Waldschmidt \cite{W2000}, Ch. 3.2.

The next lemma is due to Matveev \cite{M} and is presently the strongest general lower bound for linear forms in logarithms of algebraic numbers. He published originally its additive version, but for our purpose this multiplicative version is more beneficial.

Under the above mentioned circumstances, Matveev \cite{M} proved the additive variant of the following theorem.

\begin{Le}\label{Matveev}
Let $\mathbb{K}$ be an algebraic number field of degree $d_{\mathbb{K}}$ and let $\eta_1, \ldots, \eta_t \in \mathbb{K} \setminus \{0\}$, and $e_1, \ldots, e_t$ be nonzero integers. Put
$$
E =\max\{|e_1|, \ldots, |e_t|, 3\}\qquad{\rm and}\qquad
\Gamma = \prod_{i=1}^t \eta_i^{e_i}-1.
$$
Let $F_1, \ldots, F_t$ be such that
$$
F_j \geq \max \{d_{\mathbb{K}}h(\eta_j), |\log \eta_j|, 0.16\}, \quad {\text{\rm for}}\quad j=1,\ldots t.
$$
 If $\Gamma \neq 0$, then
\begin{equation*}
\label{ineq:matveev} \log |\Gamma| > -3\cdot
30^{t+4}(t+1)^{5.5}d_{\mathbb{K}}^2(1+\log d_{\mathbb{K}})(1+\log tE)F_1F_2\cdots F_t.
\end{equation*}
\end{Le}

The following technical lemma has many variants, its present form appeared in \cite{SL}.

\begin{Le} \label{l:n_logn}
If $s\ge 1, T\ge (4s^2)^s$, and $x/(\log x)^s < T$, then
$$
x < 2^s T (\log T)^s.
$$
\end{Le}

We also need the following simple property of the logarithm function, for its proof see de Weger \cite{deWeger}.

\begin{Le} \label{l:deWeger}
Let $a\in \mathbb{R}$. If $a<1$ and $|x|<a$ then
$$
|\log(1+x)| < \frac{-\log(1-a)}{a} |x|,
$$
and
$$
|x| < \frac{a}{1-e^{-a}}|e^x-1|.
$$
\end{Le}

\section{Proof of Theorem \ref{t:ABC_gen}}

If $u\ge u_1$ then, by  Lemma \ref{l:roots} (ii), $\alpha_{q+1}(u)$ is the dominating root of $P_u(X)$ and $1/\alpha_{0}(u)$ is the dominating root of $X^{q+2}P_u(1/X)$, hence, by the Th\'eor\'eme of Mignotte \cite{Mignotte1}, the equation \eqref{e:ABC} has for any fixed $u\ge u_1$ only finitely many effectively computable solutions. This happens for all $u_1\le u \le \max\{u_1,\ldots,u_5\}$, i.e. we may assume that the assertions of the Lemmas  \ref{l:roots}, \ref{l:a0:a_q+2}, \ref{l:Dq+1/D0v2}, \ref{l:novekedes} and of the Corollaries \ref{c:alfa2szorzat}, \ref{c:monoton} hold, moreover we have to deal only with those $u$ which are larger than $\max\{u_1,\ldots,u_5\}$.

Equation \eqref{e:ABC} compress three equations depending on the signs of $n$ and $m$. In the cases $n,m\ge 0$ and $n,m< 0$ it is enough to study the growth of $(|A_n|)$.  The situation is much more complex, if $n$ and $m$ have different signs. Because \eqref{e:ABC} is symmetric in $n$ and $m$ we may restrict ourself to $n\ge 0$ and $m<0$. We assume in the sequel that $(n,m)\in \mathbb{Z}, n\not=m$ is a solution of \eqref{e:ABC} and split the proof of Theorem \ref{t:ABC_gen} into two cases.

\medskip

{\bf Case I, $n,m\ge 0$ and  $n,m < 0$.} Assume that $|A_n(u)| = |A_m(u)|$ for some $u>u_2$. By Corollary \ref{c:monoton} the sequences $(|A_n(u)|)_{n\ge n_3}$ and $(|A_n(u)|)_{n\le -n_3}$ are strictly increasing hence \eqref{e:ABC} is not possible if $u> u_2$ and $n,m\ge n_3$ as well as if $-n,-m\ge n_3$. 

\medskip

{\bf Case II, $n\ge 0, m<0$.} The proof of this case will be divided into several steps. To simplify the notation we will indicate the dependence on $u$ only if it is essential. We start with general remarks. 
 
The inequalities \eqref{e:pozitiv} and \eqref{e:negativ} imply 
$$
-c_8 L F^n < |A_{n}|-|a_{q+1}\alpha_{q+1}^{n}| < c_8 L F^n
$$
and 
$$
-c_8 L F^{-m} < |A_{m}|-|a_{0}\alpha_{0}^{m}| < c_8 L F^{-m},
$$
where $2F<\alpha_{q+1}, |\alpha_0|^{-1}$. By our assumption \eqref{e:ABC} holds for the pair $(n,m)\in \mathbb{Z}, \; n\ge 0> m$, hence
\begin{eqnarray}
-c_8L(F^n + F^{-m}) &<& (|A_{m}|-|a_{0}\alpha_{0}^{m}|) - (|A_{n}|-|a_{q+1}\alpha_{q+1}^{n}|) \nonumber\\
 &=& |a_{q+1}\alpha_{q+1}^{n}| - |a_{0}\alpha_{0}^{m}| \label{e:also-felso}\\ 
 &<& c_8L(F^n + F^{-m})\nonumber.
\end{eqnarray}

There is a constant $n_4$ such that $|a_{q+1}|\alpha_{q+1}^{n} > 2c_8LF^n$, hence
$$
|a_{q+1}|\alpha_{q+1}^{n} < |A_n| < 3 |a_{q+1}|\alpha_{q+1}^{n}
$$
hold for all $n\ge n_4$. Similarly, if $-m\ge n_4$ than 
$$
|a_{0}\alpha_{0}|^{m} < |A_m| < 3 |a_{0}\alpha_{0}|^{m}
$$

Thus if \eqref{e:ABC} holds, and $n,-m\ge n_4$ then 
\begin{eqnarray*}
|a_{q+1}|\alpha_{q+1}^{n} &<& |A_n| = |A_m|< 3 |a_{0}\alpha_{0}|^{m}\\
|a_{0}\alpha_{0}|^{m} &<& |A_n| = |A_m|< 3 |a_{q+1}|\alpha_{q+1}^{n},
\end{eqnarray*}
which yield
\begin{equation} \label{e:alsofelsonm}
-\log|a_0/(3a_{q+1})|<  n\log \alpha_{q+1} - m \log |\alpha_0| < \log (3|a_0/a_{q+1}|).
\end{equation}

\subsection{Lower bound for $n$ and $-m$.}

So far we did not need any restrictions on $Q(0)$, apart that it is a non-zero integer. From here on we assume that $Q(0)=\pm 1$, hence $|Q(0)\alpha_0|\alpha_{q+1}= |\alpha_0|\alpha_{q+1}$. By Corollary \ref{c:alfa2szorzat} $|\alpha_0| \alpha_{q+1}$ is either larger or less than one. We treat the two cases separately.

{\bf We assume first $|\alpha_0| \alpha_{q+1}> 1$}. Inequality \eqref{e:alsofelsonm} implies 
$$
n < -m \frac{-\log|\alpha_0|}{\log \alpha_{q+1}} + \frac{\log 3|a_0/a_{q+1}|}{\log \alpha_{q+1}} < -m +  c_9,
$$
where 
$$
c_9 = (q +1)^2 + \frac{\log(3\cdot (c_4L)^{q+1})}{\log\alpha_{q+1}}.
$$
This is true because we can estimate $3|a_0/a_{q+1}| = 3 |D_0/D_{q+1}|$ by using Lemma \ref{l:Dbounds}.

The relation between $n$ and $-m$ yields $F^n<F^{-m+c_9}$, which allows us to  rewrite inequality \eqref{e:also-felso} as  

\begin{equation} \label{e:alsofelsonm1}
||a_{q+1}|\alpha_{q+1}^{n} - |a_{0}\alpha_{0}^{m}|| < c_{10}F^{-m}, \; c_{10}= c_8L(F^{c_9}+1).
\end{equation}
After division by $|a_{0}\alpha_{0}^{m}|$ we get
$$
\left| \left|\Phi(u) \right| (|\alpha_0|\alpha_{q+1})^{n-\kappa}|\alpha_{0}|^{-m-n+\kappa} - 1 \right| < \frac{c_{10}}{a_0} (F|\alpha_0|)^{-m} 
$$
with 
$$
\Phi(u) = u^{\kappa}\frac{a_{q+1}}{a_0}=u^{\kappa}\frac{D_{q+1}(u)}{D_{0}(u)}
$$
and the $\kappa\in \mathbb{Q}$ is defined in Lemma \ref{l:Dq+1/D0v2}. By Lemma \ref{l:Dbounds}
$$
|a_0| = \left|\frac{D_0}{D}\right| > c_{11} \alpha_{q+1}^{-q(q+1)}
$$
with an appropriate constant $c_{11}$. Hence
$$
\left| \left|\Phi(u) \right| (|\alpha_0|\alpha_{q+1})^{n-\kappa}|\alpha_{0}|^{-m-n+\kappa} - 1 \right| < c_{12} (F|\alpha_0|)^{-m},
$$
where $c_{12} = c_{10}\alpha_{q+1}^{q(q+1)}/ c_{11}$.

There is a constant $u_6$ such that if $u\ge u_6$ then $|\alpha_0(u)|^{-1}> c_{12} F^2$ holds. Then 
$$
|\alpha_0(u)|^{-|m|/2}> c_{12}^{|m|/2} F^{|m|}> c_{12} F^{|m|}.
$$ 
Multiplying both sides by $|\alpha_0(u)|^{|m|/2}<1$ we get 
$$
c_{12} (F|\alpha_0|)^{-m} < |\alpha_0|^{-m/2}.
$$
Thus if $u\ge u_6$ and $n\ge n_4$ then 
\begin{equation} \label{e:Baker1}
\left| \left|\Phi(u) \right| (|\alpha_0|\alpha_{q+1})^{n-\kappa}|\alpha_{0}|^{-m-n+\kappa} - 1 \right| < |\alpha_0|^{-m/2}.
\end{equation}

The number on the right hand side is obviously less than $1/2$, hence
\begin{equation} \label{e:Baker2}
\left| \log \left|\Phi(u) \right| + (n-\kappa)\log(|\alpha_0|\alpha_{q+1}) + (-m-n+\kappa) \log |\alpha_{0}| \right| < 2 |\alpha_0|^{-m/2}
\end{equation}
by Lemma \ref{l:deWeger}.

\begin{Le} \label{l:n+m}
Assume that $(n,m.u)\in \mathbb{Z}^3$ is a solution of \eqref{e:ABC} such that $n\ge 0, m<0$. Denote $\kappa_1,\kappa_2,\varphi$ the constants introduced in Corollary \ref{c:alfa2szorzat} and Lemma \ref{l:Dq+1/D0v2}. If $|\alpha_0(u)|\alpha_{q+1}(u)>1$ and $u\ge u_2$ then $n+m-\kappa\le 0$, and equality may hold only if $|\Phi(u)| <1$ and $\kappa_1>\kappa_2$ or $|\Phi(u)| <1$ and $\kappa_1 \le \kappa_2$ and $|\varphi| < 1$.  Moreover, if $u>\max\{u_1,\ldots,u_{10}\}$ then there exist effectively computable constants $g,f>0$ such that 
$$
-m \ge n-\kappa > g u^f.
$$
\end{Le}

\begin{proof}

In the proof of the lemma we separate some cases, which need different treatments.

If $n+m-\kappa> 0$ then 
$$
\left|\Phi(u) \right| (|\alpha_0|\alpha_{q+1})^{n-\kappa}|\alpha_{0}|^{-m-n+\kappa} \ge \frac{|\varphi|}{2}\left( 1+ \frac{1}{2u^{\kappa}}\right) \frac{u}{2} > 1+ \frac{1}{2u^{\kappa}}
$$
because $|\alpha_0|^{-m-n+\kappa}\ge 1/|\alpha_0|> u/2$, provided $u\ge u_2$. This contradicts \eqref{e:Baker1}, and the proof of the first statement is complete.

Assume that $n+m-\kappa< 0$. Rearranging \eqref{e:Baker2} we obtain
$$
(n-\kappa)\log(|\alpha_0|\alpha_{q+1})  > (n+m-\kappa) \log |\alpha_{0}| - \log \left|\Phi(u) \right| - 2 |\alpha_0|^{-m/2}> \left|\log \frac{ |\alpha_0|}{e(|\varphi|+1)} \right|
$$
Here we used that $(n+m-\kappa) \log |\alpha_{0}|\ge |\log |\alpha_{0}||,$ $2 |\alpha_0|^{-m/2}<1$ and $|\Phi(u)|< |\varphi|+1$ by Lemmata \ref{l:roots} and \ref{l:Dq+1/D0v2}. Corollary \ref{c:alfa2szorzat} yields $0< |\alpha_0| \alpha_{q+1}-1 < \frac{2\pi}{u^{\kappa_1}}<\frac12$ as well, but then 
$$
\log (|\alpha_0|\alpha_{q+1})< \frac{4\pi}{u^{\kappa_1}}
$$
by Lemma \ref{l:deWeger}. Hence
$$
-m > n -\kappa >\frac{u^{\kappa_1}}{4\pi} \left| \log \frac{ |\alpha_0|}{e(\varphi+1)}\right| > \frac{u^{\kappa_1}}{4\pi} \log \frac{u}{2e(|\varphi|+1)},
$$
and the stated inequality holds with $f=\kappa_1, g=(\log 2)/4\pi$, when $u>2e(|\varphi|+1)$.
\medskip

It remains to treat the case when $n+m-\kappa =0$. We distinguish some cases.

{\bf Case 1, $|\Phi(u)|\ge 1$}. As $n+m-\kappa = 0$ we have 
$$
\left|\Phi(u) \right| (|\alpha_0|\alpha_{q+1})^{n-\kappa}|\alpha_{0}|^{-m-n+\kappa} \ge |\alpha_0|\alpha_{q+1}> 1+ \frac{1}{2u^{\kappa_1}}
$$
by Corollary \ref{c:alfa2szorzat}, but this contradicts \eqref{e:Baker1} whenever $-m = n-\kappa\ge 2\kappa_1$. 

{\bf Case 2}, $|\Phi(u)| < 1$. The inequality $|\Phi(u)-\varphi|<2\mu/u^{\kappa_2}$ of Lemma \ref{l:Dq+1/D0v2} implies   
$$
|\Phi(u)| > |\phi| - ||\Phi(u)|-|\phi|| > |\varphi|-2\mu/u^{\kappa_2} >|\varphi|/2,
$$
whenever $u>u_5$.

\medskip

{\bf Subcase 2a, $\kappa_1 \le \kappa_2$}. 

If $|\varphi|=1$ then we have  
$$
|\Phi(u)|> 1- \frac{\mu}{u^{\kappa_2}} \ge 1-\frac{\mu}{u^{\kappa_1}}.
$$ 

Corollary \ref{c:alfa2szorzat} and a well known property of the exponential function imply
$$
(|\alpha_0|\alpha_{q+1})^{n-\kappa}> \left(1+ \frac{\pi}{2u^{\kappa_1}}\right)^{n-\kappa} > 1+ \frac{\pi(n-\kappa)}{2u^{\kappa_1}}.
$$
Using this inequality and $n+m-\kappa=0$ we obtain
\begin{eqnarray*}
\left|\Phi(u) \right| (|\alpha_0|\alpha_{q+1})^{n-\kappa}|\alpha_{0}|^{-m-n+\kappa} &>& \left(1-\frac{\mu}{u^{\kappa_1}} \right) \left( 1+ \frac{\pi(n-\kappa)}{2u^{\kappa_1}}\right)\\
&\ge& 1 + \frac{\pi(n-\kappa)- 2\mu}{2u^{\kappa_1}}-\frac{\pi \mu (n-\kappa)}{2u^{2\kappa_1}}\\
&\ge& 1 + \frac{\pi(n-\kappa)}{4u^{\kappa_1}}-\frac{\pi \mu (n-\kappa)}{2u^{2\kappa_1}}\\
&=& 1 + \frac{\pi(n-\kappa)}{4u^{\kappa_1}}\left(1- \frac{2\mu }{u^{\kappa_1}}\right)\\
&>&  1 + \frac{\pi(n-\kappa)}{8u^{\kappa_1}} .
\end{eqnarray*}
The third inequality holds for all $n>\frac{4\pi}{\mu}+\kappa = n_5$, while the last for all $u$ with $u>(4\mu)^{1/\kappa_1} =u_7$. The final inequality contradicts \eqref{e:Baker1}, i.e $n+m-\kappa<0$ as in Case 1. 
\medskip

If $|\varphi|<1$, then Lemma \ref{l:Dq+1/D0v2} yields 
$$
|\Phi(u)| < |\varphi| + \frac{2\mu}{u^{\kappa_2}}< 1-\delta,
$$
with $\delta = 1-|\varphi|-\frac{2\mu}{u^{\kappa_2}}>0$, provided $u > (2\pi /(|\varphi| + \delta-1))^{1/\kappa_2}=u_8$.

Corollary \ref{c:alfa2szorzat} implies $|\alpha_0|\alpha_{q+1}< 1+ 2\pi/u^{\kappa_1}$ too. Using this inequality we obtain
\begin{eqnarray*}
(|\alpha_0|\alpha_{q+1})^{n-\kappa} &<&\left( 1+\frac{2\pi}{u^{\kappa_1}}\right)^{n-\kappa}\\
&=& \left( \left( 1+\frac{2\pi}{u^{\kappa_1}}\right)^{u^{\kappa_1}/(2\pi)}\right)^{2\pi(n-\kappa)/u^{\kappa_1}}\\
&<& e^{2\pi(n-\kappa)/u^{\kappa_1}} < 1+ \frac{4\pi(n-\kappa)}{u^{\kappa_1}}.
\end{eqnarray*}
The last inequality holds whenever
$$
\frac{2\pi(n-\kappa)}{u^{\kappa_1}}<\frac{2\pi n}{u^{\kappa_1}} <1/2,
$$ 
i.e. if $n \le \frac{u^{\kappa_1}}{4\pi}$. Hence our assertion holds with $g=1/4\pi, f=\kappa_1$.

Assuming $n+m-\kappa=0$ and $n < \min \{\frac{u^{\kappa_1}}{4\pi},\kappa + \frac{\delta}{8\pi} u^{\kappa_1}\}$ and combining recent inequalities we get
\begin{eqnarray*}
|\Phi(u)| (|\alpha_0|\alpha_{q+1})^{n-\kappa} &<& \left(1 - \delta
\right) \left(  1+ \frac{4\pi (n-\kappa)}{u^{\kappa_1}}\right) \\
&<& 1-\delta + \frac{4\pi (n-\kappa)}{u^{\kappa_1}}\\
&<& 1 - \frac{\delta}{2}.
\end{eqnarray*}
This inequality contradicts \eqref{e:Baker1}, i.e. $n+m-\kappa=0$ is only possible if $n \ge \min \{\frac{u^{\kappa_1}}{4\pi},\kappa + \frac{\delta}{8\pi} u^{\kappa_1}\} \ge \frac{\delta}{8\pi} u^{\kappa_1}$. Hence, this time our assertion holds with $g=\delta=8\pi, f=\kappa_1$.

\medskip
{\bf Subcase 2b, $\kappa_1 > \kappa_2$}. 
Another consequence of the inequality $|\Phi(u)-\varphi|<2\mu/u^{\kappa_2}$ of Lemma \ref{l:Dq+1/D0v2} is that $\Phi(u)$ and $\varphi$ have the same sign for all large enough $u$. Indeed, otherwise $|\varphi|<|\Phi(u)-\varphi|<2\mu/u^{\kappa_2}$, which is absurd whenever $u> (2\mu /|\varphi|)^{1/\kappa_2}=u_9$. 

If $|\varphi|=1$ then $|\Phi(u)|<|\varphi|$ and the lower bound $|\Phi(u)-\varphi|>\mu/2u^{\kappa_2}$ of Lemma \ref{l:Dq+1/D0v2} yields $|\Phi(u)| < |\varphi| - \mu/2u^{\kappa_2} = 1- \mu/2u^{\kappa_2}$. 

If $|\varphi|<1$ then
$$
|\Phi(u)| < |\varphi| + \mu/2u^{\kappa_2} < 1- \mu/2u^{\kappa_2},
$$
but the second inequality holds only if $u> (2(1-|\varphi|)/\mu)^{1/\kappa_2} = u_{10}$. In summary we have for $0<\varphi\le 1$ uniform upper bound for $|\Phi(u)|$. Combining that with the upper bound for $(|\alpha_0|\alpha_{q+1})^{n-\kappa}$ proved in Subcase 2a, $\varphi<1$, further assuming $n+m-\kappa=0$ and 
$$
n < \min \left\{\frac{u^{\kappa_1}}{4\pi}, \kappa+ \frac{\mu}{8\pi}u^{\kappa_1-\kappa_2}\right\}
$$
we get
\begin{eqnarray*}
|\Phi(u)| (|\alpha_0|\alpha_{q+1})^{n-\kappa} &<& \left(1 - \frac{\mu}{2u^{\kappa_2}}\right) \left(  1+ \frac{4\pi (n-\kappa)}{u^{\kappa_1}}\right) \\
&<& \left(1 - \frac{\mu}{2u^{\kappa_2}}\right) \left(1 + \frac{\mu}{2u^{\kappa_2}}\right)\\
&=& 1 - \frac{\mu^2}{4u^{2\kappa_2}}. 
\end{eqnarray*}

Thus  
$$
\left| |\Phi(u)| (|\alpha_0|\alpha_{q+1})^{n-\kappa} -1 \right| = 1- |\Phi(u)| (|\alpha_0|\alpha_{q+1})^{n-\kappa}  > \frac{\mu^2}{4u^{2\kappa_2}},
$$
which contradicts again \eqref{e:Baker1}. To avoid contradiction $n$ must be large, more precisely
$$
n \ge \min \left\{\frac{u^{\kappa_1}}{4\pi}, \kappa+ \frac{\mu}{8\pi}u^{\kappa_1-\kappa_2}\right\} \ge \frac{\mu}{8\pi}u^{\kappa_1-\kappa_2}.
$$
Our inequality holds with $g=\mu/8\pi, f=\kappa_1-\kappa_2>0$.

All assertions of the lemma are proved.
\end{proof}

We are left to deal with $|\alpha_0(u)|\alpha_{q+1}(u)<1$.

\begin{Le} \label{l:n+m_ketto}
Assume that $(n,m.u)\in \mathbb{Z}^3$ is a solution of \eqref{e:ABC} such that $n\ge 0, m<0$. Denote $\kappa_1,\kappa_2,\varphi$ the constants introduced in Corollary \ref{c:alfa2szorzat} and Lemma \ref{l:Dq+1/D0v2}. If $|\alpha_0(u)|\alpha_{q+1}(u)<1$ and $u\ge u_2$ then $n+m-\kappa\ge 0$, and equality may hold only if $|\Phi(u)| >1$ and $\kappa_1>\kappa_2$ or $|\Phi(u)| >1$ and $\kappa_1 \le \kappa_2$ and $|\varphi| < 1$.  Moreover, if $u>\max\{u_1,\ldots,u_{10}\}$ then there exist effectively computable constants $g,f>0$ such that 
$$
n \ge -m+\kappa > g u^f.
$$

\end{Le}

\begin{proof}
If $|\alpha_0(u)|\alpha_{q+1}(u)<1$ then inequality \eqref{e:alsofelsonm} implies 
$$
-m < n \frac{\log \alpha_{q+1}}{-\log|\alpha_0| } + \frac{\log |a_0/(3a_{q+1})|}{-\log \alpha_{0}} < n +  c_9,
$$
with the constant $c_9$ defined before Lemma \ref{l:n+m}. Following the argumentation there we now derive from inequality \eqref{e:also-felso} with division by $|a_{q+1}|\alpha_{q+1}^{n}|$ and performing some manipulation
\begin{equation} \label{e:Baker3}
\left||\Phi(u)|^{-1}(|\alpha_0(u)|\alpha_{q+1}(u))^{^{-n+\kappa} } \alpha_{q+1}(u)^{n+m -\kappa} -1   \right| < \alpha_{q+1}(u)^{-n/2},
\end{equation}
which yields 
\begin{equation} \label{e:Baker4}
\left| \log \left|\Phi(u) \right| + (n-\kappa)\log(|\alpha_0|\alpha_{q+1}) + (-m-n+\kappa) \log |\alpha_{q+1}| \right| < 2 \alpha_{q+1}^{-n/2}.
\end{equation}

From here on one can follow the same argumentation as used in the proof of Lemma \ref{l:n+m} taking into account that $|\alpha_0|^{-1}$ and $\alpha_{q+1}$ have the same order of magnitude. As an example we prove that $n+m-\kappa\le 0$. 
Assume in the contrary that $n+m-\kappa> 0$. If $n>\kappa$ then 
$$
\left|\Phi(u) \right|^{-1} (|\alpha_0|\alpha_{q+1})^{-n+\kappa}|\alpha_{q+1}|^{n + m-\kappa} \ge \frac{1}{2|\varphi|}\left( 1+ \frac{1}{2u^{\kappa}}\right) \frac{u}{2} > 1+ \frac{1}{2u^{\kappa}}
$$
because $|\Phi(u)|< |\varphi| + 2|\mu|/u^{\kappa_2} < 2\varphi$ and $\alpha_{q+1}^{n+m-\kappa}\ge \alpha_{q+1}> u/2$, provided $u$ is large enough. This contradicts \eqref{e:Baker1}. 
\end{proof}
\medskip

\subsection{Final step}

We proved in the last subsection that if $n,m,u$ is a solution triplet of \eqref{e:ABC}  such that $u\ge \max\{u_1,\ldots,u_{10}\},n\ge 0,m<0$ then $n,|m| > gu^{f}$, with some effective constants $g,f>0$. We derive effective upper bound for $\max\{n,|m|\}$ by applying Matveev's theorem, Lemma \ref{Matveev}, to the inequalities  \eqref{e:Baker1} and \eqref{e:Baker3}. To this end we write \eqref{e:Baker1} in the form 
\begin{equation} \label{e:upperbound}
\left| \Gamma(u) -1\right|< |\alpha_0(u)|^{-m/2}
\end{equation}
with
$$
\Gamma(u) = \left|\frac{a_{q+1}(u)}{a_0(u)}\right| |\alpha_0(u)|^{m} \alpha_{q+1}(u)^n.
$$
With this notation \eqref{e:Baker3} becomes 
$$
\left| \Gamma(u)^{-1} -1 \right| < \alpha_{q+1}(u)^{-n/2},
$$
but then $|\Gamma(u)^{-1}|>1/2$, i.e. after division by $|\Gamma(u)^{-1}|$
\begin{equation} \label{e:upperbound2}
\left| \Gamma(u) -1 \right| < \frac{\alpha_{q+1}(u)^{-n/2}}{|\Gamma(u)|^{-1}}< 2  \alpha_{q+1}(u)^{-n/2},
\end{equation}

The left hand side of this inequality is the same as of \eqref{e:upperbound}, and there is slight difference in the right hand side. Hence the induced lower bound by Matveev's theorem is the same for \eqref{e:Baker1} and \eqref{e:Baker3}.

Our last lemma justify the applicability of Matveev's theorem.

\begin{Le} \label{l:Gamma} There is a constant $n_5$ depending only on $Q(X)$ and $L$ such that if $u\ge u_2$ then 
$$
\Gamma(u) = \left|\frac{a_{q+1}(u)}{a_0(u)}\right| |\alpha_0(u)|^{m} \alpha_{q+1}(u)^n \not= 1,
$$ 
provided $\max\{|n|,|m|\}\ge n_5$.
\end{Le}

We postpone the rather technical proof of this lemma to the end of the paper.

If $\max\{|n|,|m|\} <  n_5$ then the proof is finished. Otherwise, by Lemma \ref{l:Gamma}, Matveev's theorem is applicable to obtain lower bounds for the left hand side of the inequalities \eqref{e:Baker1} and \eqref{e:Baker3}. We have actually $\mathbb{K}=\mathbb{Q}(\alpha_0(u),\ldots,\alpha_{q+1}(u))$, hence $d_{\mathbb{K}}\le (q+2)!$. Further $t=3$, 
$$
\eta_1 = \left|\frac{a_{q+1}(u)}{a_0(u)}\right|, \; \eta_2= |\alpha_0(u)|, \; \eta_3= \alpha_{q+1}(u),
$$
and
$e_1=1,e_2=m,e_3=n$ and $E=\max\{n,|m|\}$. If $|\log(|\alpha_0(u)|)| \log(\alpha_{q+1}(u))>1$ then $n<-m+c_9 \le 2|m|$, i.e. $E\le -2m$, provided $-m\ge c_9$. If, however,  $\log(|\alpha_0(u)|) \log(\alpha_{q+1}(u)|)<1$ then $-m<n+c_9 \le 2n$, i.e. $E\le 2n$, provided $n\ge c_9$.

Plainly $h(\eta_2)=h(\eta_3)$, and
\begin{eqnarray*}
h(\eta_3) &=& h(\alpha_{q+1}(u))= \frac{1}{q+2} \left(\log(\alpha_{q+1}(u))- \log(|\alpha_0(u)|)+ \sum_{j=1}^q \max\{1,\log|\alpha_j(u)|\}\right)\\
&\le& \log(\alpha_{q+1}(u)).
\end{eqnarray*}
This is true because the sum in the large bracket is at most $3\log \alpha_{q+1}(u)$ and $q\ge 1$.

To estimate $h(\eta_1)$ we apply several times the inequalities of Lemma \ref{l:Weil height}. Let us start with
$$
h(\eta_1) = h\left(\frac{a_{q+1}(u)}{a_0(u)}\right) = h\left(\frac{D_{q+1}(u)}{D_0(u)}\right) \le h(D_{q+1}(u)) + h(D_{0}(u)).
$$
Notice that $D_{q+1}(u)$ and $D_{0}(u)$ are algebraic integers because both are the sum of $(q+2)!$ products of algebraic integers. For example $D_{q+1}(u)$ is the sum of terms of form 
$$
\alpha_0(u)^{i_0}\cdots\alpha_q(u)^{i_q}\cdot A_{i_{q+1}},
$$
where $(i_0,\ldots,i_{q+1})$ denotes a permutation of $\{0,\ldots,q+1\}$. The numbers $\alpha_j(u), j=0,\ldots,q+1$ are the zeroes of the irreducible polynomial $P_u(X)$, hence 
$$
h(\alpha_j(u))=h(\alpha_{q+1}(u))\le \log(\alpha_{q+1}(u)),\quad j=0,\ldots,q+1.
$$
Further $h(A_j)\le \log L, j=0,\ldots,q+1$. Hence
\begin{eqnarray*}
h\left(\alpha_0(u)^{i_0}\cdots\alpha_q(u)^{i_q}\cdot A_{i_{q+1}}\right) &\le & i_0 h(\alpha_0(u))+\ldots+ i_q h(\alpha_q(u)) + \log L\\
&\le& (q+1)^2 h(\alpha_{q+1}(u)) + \log L\\
&\le& c_{13}\log(\alpha_{q+1}(u)).
\end{eqnarray*}
We denote the expression of the right hand side by $H$.

Set $k = \lceil \log_2(q+2)!\rceil$, where $\log_2()$ denotes the logarithm function to base $2$. Consider the binary tree of depth $k$ with leaves the summands of $D_{q+1}(u)$. The nodes of each level represent the appropriately signed sum of its two children, hence the root is $D_{q+1}(u)$. We estimate the Weil height of the nodes, whose level is counted top down. At the first level, i.e. at the leaves, the height of each node is bounded by $H_1=H$. If the height of the nodes at the $j$-th, $j\ge 1$ level is bounded by $H_j$, then 
$$
H_{j+1}\le 2H_j + \log 2=4H_{j-1}+ 3\log 2=\ldots=2^{j} H_1 + 2^{j}-1 < 2^{j }(H+1).
$$
Hence
$$
h(D_{q+1}(u)) \le 2^{k-1}(H+1)\le (q+2)! (H+1) \le c_{14} \log(\alpha_{q+1}(u)).
$$
The same estimation holds for $h(D_{0})$, because  $D_{0}(u)$ and  $D_{q+1}(u)$ are algebraic conjugates, hence
$$
h(\eta_1) \le 2c_{14} \log(\alpha_{q+1}(u)).
$$

This means that the choice $F_1=F_2=F_3 = c_{14} \log(\alpha_{q+1}(u))$ is allowed, and Lemma \ref{Matveev} yields
\begin{eqnarray}
\log|\Gamma|&>&-c_{15} (1+\log (3E)) c_{14}^3 \log^{3}(\alpha_{q+1}(u)) \nonumber \\
&>&-c_{16} (\log E)\log^{3}(\alpha_{q+1}(u)). \label{e:lowerbound}
\end{eqnarray}
In $c_{15}$ we incorporated that $t=3$ and the estimations for $d_{\mathbb{K}}$, which depends only on $q$. 

If $|\log(|\alpha_{0}(u)|)| \log(\alpha_{q+1}(u))>1$ then we may set $E=2|m|$ and comparing this lower bound \eqref{e:lowerbound} with the upper bound \eqref{e:upperbound} we obtain
$$
\frac{|m|}{2} |\log(|\alpha_{0}(u)|)| < c_{16} \log(2|m|)\log^{3}(\alpha_{q+1}(u)).
$$
Corollary \ref{c:alfa2szorzat} yields $0<|\alpha_0|\alpha_{q+1}-1< \frac{2\pi}{u^{\kappa_1}}<\frac32$ as well, but then 
$$
\log (|\alpha_0|\alpha_{q+1})< \frac{4\pi}{u^{\kappa_1}}
$$
by Lemma \ref{l:deWeger}. Hence
$$
|\log(|\alpha_{0}(u)|)| = -\log(|\alpha_{0}(u)|) > \log (\alpha_{q+1}(u))
- \frac{4\pi}{u^{\kappa_1}}\ge \frac{\log (\alpha_{q+1}(u))}{2},
$$
which allow us to rewrite the last inequality including $m$ as
$$
\frac{2|m|}{\log2|m|}  < 8c_{16} \log^{2}(\alpha_{q+1}(u)).
$$
Now we apply Lemma \ref{l:n_logn} with $s=1$ and $T= 8c_{18} \log^{2}(\alpha_{q+1}(u))$ and get
$$
|m| < c_{17} \log^{3}(\alpha_{q+1}(u))
$$
after some simple and obvious computation.

On the other hand Lemma \ref{l:n+m} gives us a lower bound for $|m|$, which yields
$$
gu^{f} < |m| < c_{17} \log^{3}(\alpha_{q+1}(u)),
$$
with positive constants $g,f$. This inequality yields $u< c_{18}$, but then $|m|$ and $n$ are bounded too.

In the case $|\log(|\alpha_{0}(u)|)| \log(\alpha_{q+1}(u))>1$ we set $E=2n$ and compare the lower bound \eqref{e:lowerbound} with \eqref{e:upperbound2} we obtain
$$
\frac{|n|}{2} \log(\alpha_{q+1}(u)) +\log 2 < c_{16} \log(2n)\log^{3(q+1)}(\alpha_{q+1}(u)).
$$
From here we can derive upper bound for $n$ as we did for $|m|$ above. Hence our theorem holds with $u_0 = \max\{u_1,\ldots,u_{10}, c_{18}\}$ and
$$
n_0 = \max\{n_1\ldots,n_5, c_{17} \log^{3(q+1)}(\alpha_{q+1}(u_0))\}.
$$
$\Box$

{\it Proof of Lemma \ref{l:Gamma}.}
Assume that there is a $u>u_2$ such that $\Gamma(u)=1$, i.e.
\begin{equation} \label{e:multdep}
\frac{a_{q+1}(u)}{a_0(u)}\alpha_0(u)^{m} \alpha_{q+1}(u)^n = \pm1
\end{equation}
because all factors of $\Gamma$ are real. The polynomial $P_u(X)$ is irreducible, because $u>u_2$. Denote $\Sigma$ the Galois group of the normal closure of $\mathbb{Q}(\alpha_0(u))$. 

Assume that there is a $\sigma\in \Sigma$, which swaps $\alpha_0(u)$ and $\alpha_{q+1}(u)$. Then
\begin{eqnarray*}
\pm1 = \sigma\left( \frac{a_{q+1}(u)}{a_0(u)}\alpha_0(u)^{m} \alpha_{q+1}(u)^n  \right) &=& \sigma\left( \frac{D_{q+1}(u)}{D_0(u)}\alpha_0(u)^{m} \alpha_{q+1}(u)^n  \right)\\
&=& \frac{\sigma(D_{q+1}(u))}{\sigma(D_0(u))}\alpha_{q+1}(u)^{m} \alpha_{0}(u)^n \\
&=& \pm\frac{D_{0}(u)}{D_{q+1}(u)}\alpha_{q+1}(u)^{m} \alpha_{0}(u)^n\\
&=& \pm \frac{a_{0}(u)}{a_{q+1}(u)}\alpha_{q+1}(u)^{m} \alpha_{0}(u)^n.
\end{eqnarray*}
This equality together with \eqref{e:multdep} yields
$$
\left( \alpha_0(u) \alpha_{q+1}(u)\right)^{n+m} = \pm 1,
$$
but this is impossible as we showed in the proof of Lemma \ref{l:roots} $(v)$.

In the last equality chain we used Lemma \ref{l:a0:a_q+2} (iii), and that  $\sigma(D_{q+1}(u)) = \pm D_0(u)$ and $\sigma(D_0(u)) = \pm D_{q+1}(u)$, which are true by Lemma \ref{l:a0:a_q+2} $(iv)$. 

In the rest of the proof we assume that there is no $\sigma\in \Sigma$, which swaps $\alpha_0(u)$ and $\alpha_{q+1}(u)$.

If $n>0$ and $m\le 0$ then 
\begin{eqnarray*}
\left|\alpha_0(u)^{m} \alpha_{q+1}(u)^n\right| = \left|\frac{a_{0}(u)}{a_i(u)} \right| = \left|\frac{D_{0}(u)}{D_{q+1}(u)} \right| < (c_4L)^{q+1} \alpha_{q+1}(u)^{(q+1)^2}
\end{eqnarray*}
by Lemma \ref{l:Dbounds}, which implies an absolute bound, $n_{5,1}$, for $m, |n|$ immediately. 

If $n,m>0$ then we distinguish two cases: $0<m<n$ and $0<n<m$. In the first case pick a $\sigma\in \Sigma$, which maps $\alpha_{0}(u)$ to $\alpha_{i}(u)$ for some $0< i<q+1$, and $\alpha_{q+1}(u)$ to $\alpha_0(u)$. This is possible because $q\ge 1$. Applying $\sigma$ to \eqref{e:multdep} we get
\begin{equation} \label{e:utolso1}
\frac{a_{0}(u)}{a_i(u)}\alpha_i(u)^{m} \alpha_{0}(u)^n = \pm1,
\end{equation}
which implies 
$$
\left|\frac{1}{\alpha_{0}(u)}\right|^n \left|\frac{1}{\alpha_{i}(u)}\right|^m= \left|\frac{a_{0}(u)}{a_i(u)} \right| < (c_4L)^{q+1} \alpha_{q+1}(u)^{(q+1)^2}.
$$
As $|\alpha_{0}(u)\alpha_{q+1}(u)|<2$ by Corollary \ref{c:alfa2szorzat}, and $|\alpha_{i}(u)\alpha_{q+1}(u)|<F$ by Corollary \ref{c:F_becslés} we get
$$
\left(\frac{\alpha_{q+1}(u)^2}{2F}\right)^m < \left(\frac{\alpha_{q+1}(u)}{2}\right)^n \left(\frac{\alpha_{q+1}(u)}{F}\right)^m < \left|\frac{1}{\alpha_{0}(u)}\right|^n \left|\frac{1}{\alpha_{i}(u)}\right|^m.
$$
This inequality together with the former one yields an absolute bound, $n_{5,2}$ for $m$. Equation \eqref{e:utolso1} yields 
$$
\left|\frac{1}{\alpha_{0}(u)}\right|^n = \left|\frac{a_{0}(u)}{a_i(u)} \right| |\alpha_i(u)|^m < (c_4L)^{q+1} F^{n_{5,2}} \alpha_{q+1}(u)^{(q+1)^2},
$$
which implies the searched bound, $n_{5,3}$ for $n$ as in the former cases.

In the case $0<m<n$ one picks a $\sigma\in \Sigma$, which maps $\alpha_{0}(u)$ to $\alpha_{q+1}(u)$ and $\alpha_{q+1}(u)$ to $\alpha_i(u)$ for some $0< i<q+1$. The rest is similar as in the previous case. We leave its elaboration to the kind readers.

Finally if $n<0$ then \eqref{e:multdep} yields 
$$
\frac{a_0(u)}{a_{q+1}(u)}\alpha_0(u)^{-m} \alpha_{q+1}(u)^{-n} = \pm1.
$$
This equation has the same structure as \eqref{e:multdep}, but with a positive exponent of $\alpha_{q+1}(u)$. The same argumentation is applicable as for $n>0$. 

During the proof we distinguished several cases, but in each one we established a constant, $n_{5,x}, $ which depends only on $Q(X)$ and $L$ and is such that if $\max\{|n|,|m|\}>n_{5,x}$ then $\Gamma(u)\not= \pm1$. Taking 
$n_5= \max\{n_{5,1},n_{5,2},\ldots\}$ we finish the proof. $\Box$

\bigskip
{\bf Acknowledgement} The author thanks the anonymous referee for his/her valuable comments, especially for pointing out that the proof of $|\Gamma|\not=1$ was erroneous in the submitted version.

\end{document}